\documentclass[12]{article}
\usepackage{amssymb,amsmath,amsthm,amsfonts,amssymb}
\usepackage{hyperref}

\title{\textsf
{The Conjugacy Problem in Amalgamated Products I: 
 Regular Elements and Black Holes}}
\date{December 1,   2006}
\author{\textsf{Alexandre V. Borovik}\thanks{Partially supported
by the Royal Society Leverhulme Trust Senior Research Fellowship.}
\and \textsf{Alexei G. Myasnikov}\thanks{Supported by NSF grant
DMS-0405105,
 NSERC Discovery grant RGPIN 261898, and NSERC Canada Research Chair grant.} \and \textsf{Vladimir N.
Remeslennikov\thanks{Supported by EPSRC grant GR/R29451 and by RFFI
grant 02-01-00192.}}}

\date{}

\newcommand{\bea}{\begin{eqnarray*}}
\newcommand{\eea}{\end{eqnarray*}}
\newcommand{\bq}{\begin{quote}}
\newcommand{\eq}{\end{quote}}
\newcommand{\beq}{\begin{equation}}
\newcommand{\eeq}{\end{equation}}
\newcommand{\bi}{\begin{itemize}}
\newcommand{\ei}{\end{itemize}}

\newcommand{\bd}{\begin{description}}
\newcommand{\ed}{\end{description}}

\newtheorem{corollary}{Corollary}[section]
\newtheorem{theorem}[corollary]{Theorem}
\newtheorem{lemma}[corollary]{Lemma}

\theoremstyle{definition}
\newtheorem{definition}[corollary]{Definition}
\newtheorem{question}[corollary]{Question}
\newtheorem{example}[corollary]{Example}
\newtheorem{remark}[corollary]{Remark}

\begin{document}

\maketitle

\begin{abstract}
We discuss the time complexity of the word and conjugacy search
problems for free products $G = A \star_C B$ of groups $A$ and $B$
with amalgamation over a subgroup $C$. We stratify the set of
elements of $G$ with respect to the complexity of the word and
conjugacy problems and show that for the generic stratum  the
conjugacy search problem is decidable under some reasonable
assumptions about groups $A,B,C$.
\end{abstract}

\medskip
\begin{minipage}{5in}
\tableofcontents
\end{minipage}

\section{Introduction}

\subsection{Motivation}

This is the first paper in a series of four  written on the Word and
Conjugacy Problems in free products with amalgamation and HNN
extensions. In this introduction  we mention few results from the
other parts as well when this improves the presentation of the main
concepts.

 Free products with amalgamation and HNN extensions are among
 the most studied classical constructions in algorithmic and
combinatorial group theory. Methods developed to study the Word
and Conjugacy Problems in these groups became the classical models
much imitated  in other areas of group theory. We refer to Magnus,
Karrass, and Solitar book \cite{mks} for  amalgamated free
products techniques and to Lyndon and Schupp book \cite{LS} for
HNN extensions.

In 1971 Miller  proved  that the class of free products $A \ast_C B$
of free groups $A$ and $B$ with amalgamation over a finitely
generated subgroup $C$ contains specimens with algorithmically
undecidable conjugacy problem \cite{miller1}. This remarkable result
shows that the conjugacy problem can be surprisingly difficult even
in groups whose structure we seem to understand well. In a few years
more examples of HNN extensions with decidable word problem and
undecidable conjugacy problem followed (see the book by Bokut and
Kukin \cite{bk}). The striking undecidability results of this sort
scared away any  general research on the word and conjugacy problems
in amalgamated free products and HNN extensions. The classical tools
of amalgamated products have been abandoned and  replaced by methods
of hyperbolic groups \cite{bf,km,Mikhailovskii}, or automatic groups
\cite{bgss,Eps}, or  relatively hyperbolic groups \cite{Bum,Osin}.

In this and the subsequent paper \cite{amalgam2} we make an attempt
to rehabilitate the classical algorithmic techniques to deal with
amalgams.  Our approach treats  both decidable and undecidable cases
simultaneously. We show that, despite the common belief, the Word
and Conjugacy Problems in amalgamated free products are generically
easy
  and the classical algorithms are very fast  on ``most"
 or ``typical" inputs.  In fact, we analyze the computational complexity
 of even harder  algorithmic problems which
lately attracted much attention in cryptography (see
\cite{AAG,KL,petrides}, and surveys \cite{DHorn,S}), the so-called
{\em Search Normal Form} and {\em Search Conjugacy} Problems. Our
analysis is based on recent ideas of stratification \cite{M} and
generic complexity \cite{multiplicative,KMSS}, which we briefly
discuss below.

\subsection{Stratification of the set of inputs}

We start with a general formulation of our approach to algorithmic
problems and then specify it to algorithmic problems in groups.
We follow the book \emph{Computational Complexity} of
Papadimitriou \cite{Papa} for our conventions on computational
complexity.

Let $M$ be a set with a fixed {\em size} function $\emph{size} : M
\rightarrow \mathbb{R}_{\geqslant 0}$ and  ${\mathcal A}$ a
partial algorithm with inputs from $M$. Denote by  ${\rm Dom}
{\mathcal A} \subseteq M$ the set of inputs on which ${\mathcal
A}$ halts. For $w \in {\rm Dom} {\mathcal A}$ by $T_{\mathcal
A}(w)$ we denote the number of steps  required for the algorithm
${\mathcal A}$ to halt on the input $w$. If $f:
\mathbb{R}_{\geqslant 0} \longrightarrow \mathbb{R}_{\geqslant 0}$
is a standard complexity time bound, say $f(x) = x^n$, or $f(x) =
n^x$, $n \in \mathbb{N}$, then we say that $f(x)$ is a {\em
 worst case time upper bound} for ${\mathcal
 A}$ (with respect to the size function $\emph{size}$) if there exists a
 constant $C \in \mathbb{R}$ such that for every $w \in M$
 $$T_{\mathcal A}(w) \leqslant Cf(\emph{size}(w)) + C.$$
The set
$$M_f = \{ w \in M \mid T_{\mathcal A}(w) \leqslant f(\emph{size}(w))\}$$
is called the {\em $f$-stratum} of ${\mathcal A}$.

Assume now that the set $M$ comes equipped with a (finitely
additive) measure  $\mu$ which takes values in $[0,1]$. A subset $Q
\subseteq M$ is called {\em generic} ({\em negligible}) if $\mu(Q) =
1$ ($\mu(Q) = 0$). A bound $f$ is called a {\em generic upper bound}
for ${\mathcal A}$ if the set $M_f$ is generic with respect to
$\mu$. A generic upper bound $f$ is {\em tight} if it is a minimal
(with respect to the standard order $\preceq$ on the bounds) generic
upper bound for ${\mathcal A}$ from a fixed list of upper bounds
$U$. If not said otherwise, we always assume that $U$ consists of
polynomial bounds $x^n$ and exponential bounds $n^x, n \in
\mathbb{N}$. It may happen that an algorithm ${\mathcal A}$  does
not have a tight generic upper bound.

If $f$ is a tight generic upper bound for ${\mathcal A}$ then the
stratum $M_f$ is called a {\em generic stratum}.  Sometimes it is
difficult to determine  generic strata precisely, in which case it
is convenient to replace $M_f$ by a ``large enough" part of it. To
this end we introduce the following notion. A subset $RP \subseteq
{\rm Dom}({\mathcal A})$
 is called a  \emph{Regular Part} of $M$ relative to ${\mathcal A}$ if
$RP$ is a generic subset of $M$ such that  $RP \subseteq M_f$ for
some tight generic upper bound $f$  for ${\mathcal A}$.  One can
view  $RP$ as the  set of ``algorithmically typical" inputs for
${\mathcal A}$ with respect to $\mu$, so $RP$  describes the most
typical behavior of the algorithm on $M$. The  compliment $BH = M
\smallsetminus RP$ is called a \emph{black hole}.  Clearly, the
regular part $PR$
 and the black hole $BH$ are defined up to a negligible set.    In applications  $BH$
consists of elements $w$ in $M$ for
 which either  the algorithm ${\mathcal A}$  does not work at all,
 or  $T_{\mathcal A}(w)$ is not bounded by $f(\emph{size}(w))$,
  or for some reason it is just  not known   whether $w$ is in $M_f$ or not.
 Finally, for a bound $h \in U$ we say that the regular part $RP$ of ${\mathcal A}$  has at
 most $h$ time complexity if $M_h$ is generic. In particular, we say that $RP$ is
 polynomial time if it has at most $h$ time complexity for some polynomial $h$.

In the sequent the measure $\mu$ appears either as  the asymptotic
density function on $M$ with respect to the size function
$\emph{size}$, or the  exponential distribution on $M$ which comes
from a corresponding random walk on $M$ (we refer to
\cite{multiplicative,BMS} for details). To explain this we need a
few definitions. Let $$M=\bigcup_{i=0}^\infty M_i$$ be a partition
of $M$ with respect to
 the given size function $\emph{size}:M \rightarrow \mathbb{R}$,
 thus
$$M_i = \{w \in M \mid \emph{size}(w) = i\}.$$
In this case for a subset $Q$ of $M$ the fraction
 $$\frac{\mu(Q\cap S_i)}{\mu(M_i)}$$
can be viewed as the probability of an element of $M$ of size $i$ to
be in $Q$. The limit (if it exists)
$$\rho(Q) = \lim_{i\rightarrow \infty} \frac{\mu(Q\cap M_i)}{\mu(M_i)}$$
is called the \emph{asymptotic density} of $Q$. The set $Q$ is
\emph{generic} ({\em negligible})  with respect to $\rho$ if $
\rho(Q) = 1$ ($\rho(Q) = 0$), and $Q$ is \emph{strongly generic} if
the convergence
$$\frac{\mu(Q\cap S_i)}{\mu(S_i)} \rightarrow 1$$
is exponentially fast when $i\rightarrow \infty$.

It is not hard to see that the union ${\mathcal F}$ of all generic
and negligible subsets of $M$ is an algebra of subsets of $M$ and
the asymptotic density $\rho$ is a measure on the space
$(M,{\mathcal F})$.

\subsection{Search problems in groups}
 \label{se:alg-prob}

 The Word and Conjugacy Problems are two classical algorithmic problems
 introduced by M.~Dehn in 1912.  Since then much of the research in combinatorial group
 theory was related to these problems. We
refer to surveys  \cite{AdD,miller1,miller2,RR} on algorithmic
problems in groups.

Let $H$ be a fixed group given by a finite presentation $H = \langle
X ; R \rangle$, and $M(X) = (X^{\pm 1})^*$  a free monoid over the
alphabet $X^{\pm 1}$. Sometimes, slightly abusing notations,  we
identify words in $M(X)$ with their canonical images in the free
group $F(X)$.

In general, an algorithmic problem  over $H$ can be described as a
subset $D$  of a Cartesian power $M(X)^k$ of $M(X)$. The problem is
{\em decidable} if there exists a {\em decision algorithm} $A$ which
on a given input $w \in M(X)^k$ halts and outputs ``Yes" if $w \in
D$, otherwise it outputs ``No". Notice, that, on the first glance,
decidability of the problem depends on the given presentation
$\langle X ; R \rangle$ of $H$. However, if the problems under
consideration are  "algebraic", i.e., concerning
 elements, subgroups, automorphisms, etc., then  they are
 decidable relative to one finite presentation if and only if they are
  decidable relative to any finite presentation.

 Most of the algorithmic problems for groups come in three
variations: {\em specific, uniform,} and {\em search}.

The classical formulations of the Word and Conjugacy Problems are
specific,  i.e.,  for a fixed presentation $H = \langle X ; R
\rangle$ find an algorithm (specific to this group presentation) to
solve these problems in $H$. However, this particular meaning of an
algorithmic problem has changed somewhat in recent years under
influence of practical computations with groups. Development of
software packages for computing in groups often requires
implementation of {\em uniform} decision algorithms that are able to
deal with a given algorithmic problem in various classes of
presentations of groups. This might change the flavor and
computational complexity of the problem dramatically. For example,
the uniform Word Problem  for a class of presentations $\mathcal P$
 has a pair:  a  presentation $\langle X ; R
\rangle \in \mathcal{P}$ and a word $w \in F(X)$ -   as its input,
and it requires to verify whether or not $w$ is equal to 1 in  the
group $H = \langle X ; R \rangle$. To see the difference, assume
that one solves the (specific)  Word Problem in a nilpotent group $H
= \langle X ; R \rangle$. Since $H$ has a faithful matrix
representation  $\rho_H: H \rightarrow UT_n(\mathbb{Z})$ for a
suitable $n$ one can replace the generators $x \in X$ by their
images $\rho_H(x)$ and perform matrix multiplication to check if
$\rho_H(w) = 1$ or not. In the case of the uniform variation of the
Word Problem the algorithm above would require first to find a
faithful representation $\rho_P$ for a given finite presentation $P
\in \mathcal{P}$ - not easy task in itself. Similarly, if the group
 $H = \langle X ; R \rangle$ is hyperbolic then one can use the
 famous Dehn algorithm to solve the (specific) Word Problem in $H$.
 Namely, let  $P' = \langle X ; R'
 \rangle$  be an arbitrary finite Dehn presentation of $H$
   (it always exists in this case). It is known that  $w = 1$ in $H$ if and only if the
   Dehn algorithm relative to $P'$   rewrites the word $w$ into the empty
   word. The described algorithm is very fast (linear time in the
   length of $w$, see \cite{Anshel}), but it relies on the knowledge
   of the Dehn Presentation $P'$.  If $\mathcal{P}$ is the class of
   all finite presentations of hyperbolic groups then the uniform
   version of the algorithm would require for a given presentation $P \in \mathcal{P}$
    to compute a Dehn presentation $P'$ for the group $H$ - which is
    again very demanding (there are no polynomial time algorithms known at
    the moment).

In this paper we study uniform algorithms for the Word  and
Conjugacy Problems in the class of groups which are free products
with amalgamation given by their standard presentations.

Observe, that the uniformity of the problem may appear on different
levels, not related to sets of presentation at all.
 For example, the  specific  Membership Problem
in $H$ is decidable for a given fixed finitely generated subgroup
$D$ of $H$ if there exists an algorithm which for every word $w \in
F(X)$ decides whether the element represented by $w$ in $H$ belongs
to $D$ or not.   Meanwhile, decidability of the uniform Membership
Problem for $H$ requires an algorithm which would solve the specific
Membership Problem for every finitely generated subgroup $D$ of $H$.

Finally, the {\em search} variation of an algorithmic problem  $D$
requires to decide whether a given $w$ belongs to $D$ or not, and if
it belongs, to provide   a ``reasonable proof" that $w$ is, indeed,
in $D$. For instance, the Search Word Problem  for  $H = \langle X ;
R \rangle$  usually requires to check if $w \in F(X)$ is equal to
one in $H$, and if so,  represent $w$ as a product of conjugates of
relators from $R$. The Search Word Problem is sometimes provably
harder  then the solution of the original Word Problem. Indeed, the
group $BS(1,2) = \langle a,b ; a^{-1}ba = b^2\rangle$ has a
polynomial time decision algorithm for the  Word Problem, but its
Dehn function is exponential \cite{Gersten}, so it requires, in the
worst case, at least exponential time to represent $w$ as a product
of conjugates of the relator. This new requirement for the search
decision problems to provide a ``proof", or a ``witness", of the
correct decision brings
 quite a few new algorithmic aspects, which were not studied
  in group theory. We refer to \cite{KMSS,MU} for a more
detailed discussion of the Search Problems in groups.

 Search problems could
also be uniform or specific. It is convenient to treat  uniform and
specific forms as particular cases of problems which are {\em
uniform relative  to a given class of objects} $\Phi$. More
precisely, let $D$ be an algorithmic problem on a set of inputs $I$.
We say that $D$ is decidable on a subset $\Phi \subseteq I$ if there
exists a partial algorithm $\mathcal A$ with a halting set
$Dom({\mathcal A}) \subseteq I$ that correctly solves $D$ on every
input from $Dom({\mathcal A})$ and such that $\Phi \subseteq
Dom({\mathcal A})$.

For example, the membership problem for a class of subgroups $\Phi$
of $H$ solves the specific membership problem for every subgroup $D$
from $\Phi$. If the set $\Phi$ is the whole set of elements,
finitely generated subgroups, etc., of $H$ then we omit it from the
notation. This relative approach is very natural, since there are
groups in which the uniform version of a particular algorithmic
problem is undecidable, but still there are interesting subclasses
of objects $\Phi$ for which this problem is uniformly decidable.
Moreover, even if the uniform version of the problem is decidable
the class of all objects in the question can be partition into
different subclasses with respect to different complexities  of the
decision algorithms.

Below we list some  algorithmic problems for $H$  in their uniform
relative to a subclass search variation. These algorithmic problems
involve different subsets of $H$ (subgroups, cosets, double cosets,
regular sets, etc.) given by some natural effective (constructive)
descriptions. For example, finitely generated subgroups $D$ are
given by finite generating sets (which are given as  words from
$F(X)$), cosets $wD$ are given as pairs $(D,w)$, regular sets are
given either  by finite automata or by regular expressions, etc.
Usually,  we do not specify any particular descriptions of these
subsets, unless it is required by complexity issues or by a
particular algorithm.

The Search Word Problem for finitely presented groups has several
formulations which depend on the  form of the witness. The following
is, perhaps, the most typical one.

\medskip\noindent
{\bf The Word Search Problem for a given subset of elements $\Phi$
($\mathbf{WSP_{\Phi}}$):} {\em Let $\Phi$ be a given subset of
elements from $H$ (given as words from $F(X)$). For a given  $w \in
\Phi$ decide whether $w = 1$ in $H$ or not? If $w = 1$  then find a
presentation of $w$ as a product of conjugates of relators from
$R$.}

  However, in free
products with amalgamation it is convenient to consider the
following variation of the Search Word problem.  Let ${\mathcal N}$
be a fixed set of normal forms (viewed as words in $M(X)$) of
elements from $G$, and ${\bar w}$ a representative of $w$ in
${\mathcal N}$.

\medskip\noindent
{\bf The Normal Forms Search Problem for a given subset of elements
$\Phi$ ($\mathbf{NFSP_{\Phi}}$) of $H$:} {\em Let $\Phi$ be a given
subset of elements from $H$ (given as words from $F(X)$). For a
given  $w \in \Phi$   find its normal form ${\bar w} \in {\mathcal
N}$. }

\medskip\noindent

  In free products with amalgamation and HNN extensions it is
  convenient to begin with  the {\em reduced forms} of elements
  and then specify the normal forms and cyclically reduced normal forms among them.
  One can introduce similarly  the Search Problems for reduced and cyclically reduced normal forms.
  We leave it to the reader.

\medskip\noindent
{\bf The Conjugacy Search  Problem for a given set of pairs of
elements $\Phi$ ($\mathbf{CSP_{\Phi}}$):} {\em Let $\Phi$ be a given
set of pairs of elements from $H$. For a given pair $(u, v) \in
\Phi$ determine whether $u$ is a conjugate of $v$ in $H$ or not, and
if it is then find a conjugator.}

\medskip\noindent
{\bf  The Membership Search Problem for a set of subgroups $\Phi$
($\mathbf{MSP_{\Phi}}$) of $H$:} {\em Let $\Phi$ be a set of
finitely generated subgroups of $H$. For a given $D \in \Phi$ and a
given $w \in F(X)$ determine whether  $w$ belongs to $D$ or not, and
if so, find a decomposition of $w$ as a product of the given
generators of $D$. }

\medskip\noindent
{\bf  The Conjugacy Membership Search Problem for a set of subgroups
$\Phi$ ($\mathbf{CMSP_{\Phi}}$):} {\em Let $\Phi$ be a set of
finitely generated subgroups of $H$. For a given $D \in \Phi$ and a
given $w \in F(X)$ determine whether  $w$ is a conjugate of an
element from $D$, and if so, find such an element in $D$ and a
conjugator.}

\medskip\noindent
{\bf  The Coset Representative Search Problem for a set of subgroups
$\Phi$  ($\mathbf{CRSP_{\Phi}}$):} {\em Let $\Phi$ be a set of
finitely generated subgroups of $H$.  For a given $D \in \Phi$ find
a recursive set $S$ of representatives of $D$ in $H$ and an
algorithm $A_S$ which for a given word $w \in F(X)$ finds a
representative for $Dw$ in $S$}

\medskip

Observe that to solve $\mathbf{CRSP_\Phi}$ for a given $D \in \Phi$
it suffices to find the algorithm $A_S$, since $w \in S$ if and only
if $w$ is the output of $A_S$ on the input $w$.

To formulate the next algorithmic problem we need the following
definition.  Let $M$ be a subset  of a group $H$. If $u,v \in H$
then the set $uMv$ is called a shift of $M$. For a set $\mathcal{M}$
of subgroups of $H$ denote by $\Phi({\mathcal{M}},H)$ the least set
of subsets of $H$ which contains ${\mathcal{M}}$ and is closed under
shifts and intersections.

\medskip\noindent
{\bf The Cardinality Search Problem for $\Phi(\mathcal{M},H)$
($\mathbf{CardSP_{\Phi}}$):} {\em Let ${\mathcal{M}}$ be a
collection of subsets  of $H$. Given a set $D \in
\Phi({\mathcal{M}},H)$ decide whether $D$ is empty, finite, or
infinite and, if $D$ is finite non-empty, list all elements of $D$.}

\medskip\noindent

\subsection{Results}

We show below that, under some reasonable assumptions about the
groups in the free product with amalgamation $G = A \ast_C B$, the
Normal Forms Search Problem for the classical
 normal forms in $G$ is decidable and the Search Conjugacy Problem
 is decidable for the set of {\em regular} elements in $G$.
  Moreover, we analyze the time complexity of these problems (modulo the
 corresponding algorithms in the factors $A$, $B$).

 In Section \ref{se:can} we study time complexity   of the
 standard Algorithm II  for computing the normal forms of elements
 in amalgamated free products (as, basically,  described
 in \cite{mks}).

The direct inspection of  the algorithms reveals two key issues that
determine its  time complexity.
 The  first one is related to the time complexity of the decision algorithms of
 the relevant   algorithmic problems that have to be decidable in the factors
 (the Subgroup Membership Search Problem  and
 the Subgroup Representatives Search Problem).
The second issue concerns with the following phenomena that occurs
when computing the required form (normal or reduced) for the input
$w=g_1 g_2\cdots g_n$, $g_i\in A\cup B$ in the group $A \ast_C B$.
While in progress, Algorithm  I gradually  rewrites the words in
generators of $A$ and $B$ representing the elements $g_i$, possibly
increasing their length. It happens sometimes that the accumulated
increase in length is exponential in terms of the original length.

These two issues are very different in nature, the algorithmic
difficulties of the first type come with the factors and we view
them as part of the given data. Meanwhile, the ones of the second
type are  intrinsic  to the construction itself. To deal with the
first issue we elected to specify precisely which {\em basic
algorithmic problems} are required to be decidable in the factors
and study algorithms for $A \ast_C B$ "modulo" the time complexity
of the basic problems in the factors $A, B$.  In this case  one can
view every instance of execution of a basic algorithm (a decision
algorithm corresponding to a basic problem) as one "elementary
step". It turns out that  if the complexity of the basic algorithms,
as well as the intrinsic complexity,  is known then one can estimate
the total complexity of  Algorithms I.

As an example of this kind of analysis we give the following result
for free products of free groups in Section \ref{se:comp-Alg-I}.

\newpage
{\bf Theorem 3.7} {\em
\begin{enumerate}
\item [{\rm (1)}] Let $A \ast_C B$ be a free product of free groups with
finitely generated amalgamated subgroup $C$. Then Algorithm II has
at most exponential  time complexity function (in the syllable
length of the input words);
 \item [{\rm (2)}] There are finitely generated free groups  $A$ and
$B$ and a finitely generated  subgroup $C$ in $A$ and $B$ for which
 the lower (and upper) bounds  on the time complexity of Algorithm II are exponential.

\end{enumerate}
}

\medskip

 In Section \ref{se:4} we study the time complexity of the Search
 Conjugacy Problem in $G = A \ast_C B$. More precisely,  we study
 the time complexity of  the search variation of the standard
 decision algorithm for the Conjugacy Problem
 in amalgamated free products (following the description
 in \cite{mks}). The main result of this paper shows that this
 algorithm   solves the Search
 Conjugacy Problem in $G$ for all   cyclically reduced {\em regular elements}, and
 their  conjugates, provided the relevant basic algorithmic problems  in the factors are
 decidable. We would like to emphasize that the algorithm is {\em partial}, i.e.,
 it halts and gives the  answer (always the correct one) only for inputs from a subset of $G$
 (though a very big one). This is a crucial aspect of our approach, since, as we have mentioned above,
the Conjugacy Problem  in some of these  groups is  undecidable,
hence  there are no total decision algorithms for the Conjugacy
Problem in these groups. We show in the subsequent papers that this
partial algorithm is, perhaps, as good as a total one, and may be
even better on the most typical inputs (if the total algorithm
becomes very inefficient when trying to accommodate  all non-regular
inputs).

To describe the regular elements in $G$ we
 need the following definitions. The set
 $$ N^*_G(C) = \{g \in \mid C\cap C^g \ne 1\}$$
is called the  {\em generalized normalizer} of $C$ in $G$
\cite{bmr}. Its "dual" is defined by
 $$Z_G(C) = \bigcup_{g \in N^*_G(C)\smallsetminus C} C^{g^{-1}} \cap C$$

  Now in the group $G = A \ast_C B$ the set
 $$BH = ({N^\ast}_G(C) \smallsetminus C) \cup  {Z}_G(C) $$
 is called a {\it Black Hole}, and its complement $RP = G \smallsetminus BH$ is
 called the {\it Regular Part} of ${\mathcal A}$, meanwhile the elements from $RP$
 are called {\em regular}.

The following result shows that under reasonable assumptions on the
factors the Membership Problem for the set of regular elements is
decidable in $G$.

\medskip
{\bf Theorem 4.1} {\em  Let $G = A \ast_C B$ be a free product of
finitely presented groups $A$ and $B$ amalgamated over a finitely
generated subgroup $C$. Assume also that $A$ and $B$ allow
algorithms for solving the following problems:

\begin{itemize}
 \item The Search Membership Problem for the subgroup $C$.
  \item Coset Representative Search Problem for the subgroup $C$.
   \item Cardinality Search Problem for $\Phi(Sub(C),A)$ in $A$ and for
   $\Phi(Sub(C),B)$ in $B$.
    \item The Membership Problem for $N^\ast_A(C)$, $N^\ast_B(C)$ and  $Z_A(C)$,
    $Z_B(C)$.
\end{itemize} Then there exists an algorithm to determine whether a given
element in $G$ is regular or not.}

\medskip
{\bf Theorem 4.18} {\em  Let $G = A \ast_C B$ be a free product of
finitely presented groups $A$ and $B$ amalgamated over a finitely
generated subgroup $C$. Assume also that $A$ and $B$ allow
algorithms for solving the following problems:

\begin{itemize}
  \item The Membership Search Problem for the subgroup $C$.
  \item Coset Representative Search Problem for the subgroup $C$.
  \item Cardinality Search Problem for $Sub(C)$  (see Section \ref{se:alg-prob}).
  \item Conjugacy Search Problem.
  \item Conjugacy Membership Search Problem for $C$.
\end{itemize}

Then the Conjugacy Search Problem in $G$ is decidable for cyclically
reduced regular elements from $G$ and their conjugates.}

\medskip
{\bf Corollary 4.19} {\em  Let $G = A \ast_C B$ be  a free product
of free groups $A$ and $B$ with amalgamated finitely generated
subgroup $C$. Then the Conjugacy Search Problem in $G$ is decidable
for cyclically reduced regular elements and their conjugates. }

\medskip
It is worthwhile to mention that the conjugacy problem for elements
of the syllable length greater than 1  is somewhat easier (it
requires less conditions on the factors).

\medskip
{\bf Theorem 4.15} {\em
 Let $G = A \ast_C B$ be a free product of finitely
presented groups $A$ and $B$ amalgamated over a finitely generated
subgroup $C$. Assume also that $A$ and $B$ allow algorithms for
solving the following problems:

\begin{itemize}
  \item The Membership Search Problem for the subgroup $C$.
  \item Coset Representative Search Problem for the subgroup $C$.
  \item Cardinality Search Problem for $Sub(C)$.
\end{itemize}

Then the Conjugacy Search Problem in $G$ is decidable for cyclically
reduced  regular elements $g$ of syllable length $l(g) > 1$. }

At the end of this section we briefly discuss some connections with
the subsequent papers of the series.

In the paper \cite{amalgam2} we give asymptotic estimates of the
size of the regular part $RP$ and the black hole $BH$ in free
products of  free groups with amalgamation. This enables us to show
that Algorithm II, as well as the algorithm that solves the Search
 Conjugacy Problem in such groups,  both have a
 polynomial  generic case time complexity.

 In the paper \cite{HNN1} the main results on the normal forms and the Conjugacy Search Problem
 in the free products with amalgamation are generalized to HNN extensions.
 To this end we introduce the so called {\em transfer machines}
 which allow one to "transfer" effectively the results on the classical algorithmic
 problems for free products with amalgamation to the corresponding
 HNN extensions of groups.

In the paper \cite{HNN2} we apply the general results obtained in
the papers \cite{amalgam2,HNN1} to the Miller's groups which are
particular types of HNN extensions. In particular, we show that,
despite the conjugacy problem is undecidable in these groups, there
exists an algorithm that solves the Conjugacy Search Problem in the
Miller's groups in polynomial time on "most inputs".

\section{Preliminaries}

\subsection{Amalgamated products}
\label{sec:1.1}

In this section we briefly discuss  definitions and some known facts
about free products with amalgamation. We refer to  \cite{mks} for
details.

 Let $A$, $B$, $C$ be groups and
$\phi:C \longrightarrow A\; \hbox{ and } \; \psi: C
\longrightarrow B $ be  monomorphisms. Then one can define  a
group $$ G = A *_C B,
$$ called the {\em amalgamated product of $A$ and $B$ over
$C$}  (the monomorphisms $\phi, \psi$ are suppressed from
notation). If $A$ and $B$  are given by  presentations
 $$A= \langle X \mid R_A = 1\rangle,  \ \ B=\langle Y \mid R_B =1 \rangle,$$
  and a generating set $Z$ is given for the group $C$,  then the group $G$
has presentation
\begin{equation}
\label{eq:0} G = \langle X \cup Y \mid R_A = 1,  R_B = 1,  z^\phi
= z ^\psi (z \in Z)\rangle.
\end{equation}
Notice that if the groups $A$ and $B$ are finitely presented and $C$
is finitely generated then the group $G$ is finitely presented. If
$Z = \{z_1, \ldots,z_n\}$ and we denote
$$
z_i^\phi = u_i(X), \;\; z_i^\psi = v_i(Y)
$$
then $G$  has presentation
$$
 G = \langle X \cup Y \mid R_A = 1, R_B = 1, u_1 = v_1, \ldots, u_n = v_n \rangle.
$$

The groups $A$ and $B$ a called {\em factors} of the amalgamated
product $G = A *_C B$, they are isomorphic to the subgroups in $G$
generated respectively by $X$ and   $Y$. We identify $A$ and $B$
with these subgroups via the identical maps $x \rightarrow x$, $y
\rightarrow y$ $  (x \in X,\; y \in Y)$. If we put
 $$
  V = \{u_1, \ldots, u_n\}, \;\; V = \{v_1, \ldots,v_n\}
$$
then
$$C
= \langle U \rangle  = \langle V\rangle = A\cap B \leqslant G.$$

\subsection{Normal forms of elements}
\label{sec:1.2}

Let $G = A \ast_C B$ be an amalgamated product of groups as in
(\ref{eq:0}). Denote by $S$ and $T$ fixed systems of right coset
representatives of $C$ in $A$ and $B$. Throughout this paper we assume
that the representative of $C$ is the identity element $1$.

The following notation will be in use throughout the paper. For an
element $ g \in (A \cup B) \smallsetminus C$ we define $F(g) = A$
if $g \in A$ and $F(g) = B$ if $g \in B$.

\begin{theorem} {\rm \cite[Theorem~4.1]{mks}}
An arbitrary element\/ $g$ in \/ $G = A*_C B$
 can be uniquely written in the {\em normal form}
with respect to $S$ and $T$
\begin{equation}
g = cg_1g_2\cdots g_n, \label{eq:normal}
\end{equation}
where  $c \in C$, $g_i \in T \cup S \smallsetminus \{1\}$,
 and $F(g_i) \neq F(g_{i+1})$, $i = 1,\dots, n$, $n \geqslant 0$. \label{th:1.1}
\end{theorem}

\begin{corollary} \label{cor:1.2}
Every element\/ $g \in A*_C B$ can be written in a {\em reduced
form}
\begin{equation}
\label{eq:reduced}
 g = cg_1g_2\cdots g_n
\end{equation}
where  $c \in C, g_i \in (A \cup B) \smallsetminus C$,  and
$F(g_i) \neq F(g_{i+1})$, $i = 1,\dots, n$, $n \geqslant 0$. This
form may not be unique, but the number $n$ is uniquely determined
by $g$.  Moreover, $g = 1$ if and only if\/ $n = 0$ and $c = 1$.
\end{corollary}

Let $g \in A \ast_C B$ and  $g = cg_1g_2\cdots g_n$ be a reduced
form of $g$. Then the number $n$ is  called the {\em length} of
$g$ and it is denoted by $l(g)$. Observe, that $l(g) = 0
\Longleftrightarrow g \in C.$

\begin{definition}
Let $g \in A \ast_C B$. A reduced  form  $g = cg_1g_2\cdots g_n$
is called {\em  cyclically reduced }  if one of the following
conditions is satisfied:
\begin{enumerate}
\item [{\rm (a)}] $n=0$; \item [{\rm (b)}] $n=1$ and $g$ is not a
conjugate of an element in $C$; \item [{\rm (c)}] $n \geqslant 2$
and  $F(g_1) \neq F(g_n).$
\end{enumerate}
\end{definition}
Notice that our definition of cyclically reduced forms is slightly
different from the standard one (see, for example, \cite{mks}).
Usually,  the condition (b) is not required, but the difference is
purely technical,  and it is convenient to have (b) when dealing
with conjugacy problems. Observe also, that if one of the reduced
forms of $g$ is cyclically reduced then all of them are cyclically
reduced. In this event, $g$ is called  {\em cyclically reduced}
element.

\begin{lemma}[\cite{mks}]
\label{le:cyc-red} Let $g \in  A \ast_C B$. Then $g$ is a conjugate
of  some  element $g_0$ in a cyclically reduced normal form. This
element $g_0$ is not uniquely defined, but its length $l(g_0)$ is
uniquely determined by $g$.
\end{lemma}

The normal form of $g_0$ is called {\em a cyclically reduced normal
form} of $g$.  The uniquely determined number $l(g_0)$ is called the
{\em cyclic length} of $g$ and it is denoted by $l_0(g)$. Observe
that \begin{eqnarray*} l_0(g) = 0 &\Longleftrightarrow& \ \hbox{
some  conjugate  of }\ g \ \hbox{ is
in } \ C,\\
l_0(g) = 1 &\Longleftrightarrow&  \ \hbox{ some  conjugate  of }\
g \ \hbox{ is  in } \ (A \cup B) \smallsetminus C.
\end{eqnarray*}

\subsection{The conjugacy criterion}

\begin{theorem} \cite[Theorem~4.6]{mks}
Let\/ $G = A \ast_C B$ be an amalgamated product, and  let\/ $g$
be a cyclically reduced element in $G$.

\begin{itemize} \item [(i)] If\/ $l_0(g) =0$, i.e., $g \in C$,  and $g$ is
conjugate to  an element $c \in C$ then there exists a sequence of
elements\/ $c =c_0,c_1,\dots, c_t=g$, where $c_i \in C$ and
adjacent elements $c_i$ and $c_{i+1}$, $i = 0,\dots, t-1$,
 are conjugate in $A$ or in $B$.

\item [(ii)] If\/ $l_0(g) = 1$, i.e., $g \in A \cup B \smallsetminus
C$, and\/ $g'$ is a cyclically reduced element which is  a conjugate
of $g$ in\/ $G$ then $l(g') = 1$, $F(g) = F(g')$ and  $g$ and\/ $g'$
are conjugate in $F(g)$.

\rm{(iii)} Let \/ $l_0(g) =r \geqslant 2$ and\/ $g = g_1\cdots g_r$
be  a cyclically reduced form of $g$. Assume that\/ $g$ is a
conjugate of a cyclically reduced element\/ $h = h_1\cdots h_s$ in
$G$. Then\/ $r=s$ and\/ $h$ can be obtained from $g$ by a cyclic
permutation of the elements $g_1,\dots, g_r$ followed by a
conjugation by an element from\/ $C$. \end{itemize}
\label{th:conjugacycriterion}
\end{theorem}

\subsection{Malnormal subgroups}
\label{se:maln}

 Recall, that a subgroup $H$ of a group $G$ is
called {\em malnormal } in $G$ if $H \cap H^g =1 $ for all $g \in
G \smallsetminus H$.

It follows immediately from the conjugacy criterion ( Theorem
\ref{th:conjugacycriterion})  that free factors $A$ and $B$ are
malnormal in the free product $A \ast B$. It is known that maximal
abelian subgroups (= proper centralizers) are malnormal in
torsion-free hyperbolic groups, in particular in free groups. We
refer to \cite{GKM} for results on malnormality of  maximal
abelian groups in free products with amalgamation and HNN
extensions.

\begin{definition}
Let $G$ be a group and $H$ be a subgroup of $G$. The {\em generalized
normalizer} $N^*_G(H)$ is a set of all elements $g \in G$ such that
$H\cap H^g \ne 1$.
\end{definition}

 Notice that, $N_G(H) \subseteq N^*_G(H) $, and,  in general, $N^*_G(H)$
is not a subgroup. It is obvious that  if $g \in N^*_G(H)$ then
$N^*_G(H)$ contains the whole  double coset $HgH$. A set of
representatives $\{g_i  \mid i \in I \}$ of double cosets of $H$
is called a  {\em double transversal} of $H$ in $N^*_G(H)$, in
this event
$$ N^*_G(H) = \bigcup_{i \in I} Hg_iH $$

 If $H$ is a finitely generated
subgroup of  a free group $G$  then $H$ has  a finite double
transversal in $N^*_G(H)$, moreover such a transversal can be
found algorithmically \cite{bmr}. A more convenient algorithm (in
terms of subgroup graphs) can be found in \cite{km}.

For an element $g \in G$ define
 $$Z_g(H) = \{h \in H \mid h^g \in H\} = H^{g^{-1}} \cap H$$
and put
 $$Z_G(H) = \bigcup_{g \in N^*_G(H)\smallsetminus H} Z_g(H).$$
Even though $Z_g(H)$ is a subgroup of $G$ for every $g \in G$, the
set $Z_G(H)$ may not be a subgroup. Observe, that for any $u, v
\in H$
$$Z_{ugv}(H) = {Z_g}^{u^{-1}}(H).$$
Hence if $T$ is a double transversal of $H$ in $N^*_G(H)$ then
$Z_G(H)$ is   union of conjugacy classes:
$$Z_G(H) = \bigcup_{h \in H, t \in T, t \neq 1}Z_t(H)^h.$$
In particular, if the transversal $T$ is finite then $Z_G(H)$ is
union  of  finitely many conjugacy classes of subgroups $Z_t(H)$.

\begin{definition}
Let $G$ be a group equipped with a map $L:G \rightarrow \mathbb{N}$
and $H$ be a subgroup of $G$.  For an element $g \in G$ define
$L_H(g)$ as the minimal value  of $L$ on the double coset $HgH$.
Then the {\em malnormality degree} $md(H)$ of $H$ in $G$ with
respect to $L$ is  the smallest number  $r$ such that $H \cap H^g =1
$ for all $g\in G$ with $L_H(g) \geqslant r$, if it exists,  and
$\infty$ otherwise.
\end{definition}

For example, the malnormality degree of subgroups can be  defined
in free groups, free products with amalgamation, and HNN
extensions of groups with respect to the canonical length
functions. In the sequel we always assume that for $H \leqslant A
\ast_C B$ the degree $md(H)$ is viewed with respect to the
canonical length function $l:A \ast_C B \rightarrow \mathbb{N}.$

Obviously, if a subgroup $H$ has a finite double  transversal in
$N^*_G(H)$ then $md(H)$ is finite.

\begin{lemma}
\label{le:mal-deg}
 Let\/ $G = A*_C B$ and $D \leqslant C$. Then

\begin{itemize} \item [(i)] If\/ $C$ is malnormal in $A$ and\/ $B$ then $md(D) = 1$.

\item [(ii)] If\/ $C$ is malnormal in one of the groups  $A$ and $B$
then  $md(D) \leqslant 2$. \end{itemize}
\end{lemma}

\begin{proof} Let  $g = g_1\cdots g_n$ be a reduced form of an element $g
\in G$. Suppose  $l(g) \geqslant 1$, in particular, $g_n \not \in
C$. Suppose also that $c, c' \in C$. If
$$
  g_1\cdots g_n c g_n^{-1}\cdots g_1^{-1} = c'
$$
then $g_n cg_n^{-1} \in C$.  Assume that $C$ is malnormal in both $A$
and $B$. This implies that $g_n \in C$ -contradiction.  Then $n=0$ and
therefore $md(D)=1$.

If $C$ is malnormal either in  $A$ or in $B$ then similar argument
shows that $md(D) \leqslant 2$. \end{proof}

\begin{question}
Let $G = A \ast _C B$ and $H$ be a finitely generated subgroup of
$G$. Assume that the malnormality degree $md_G(C)$ of $C$ in $G$
is finite,  and $H$ contains no elements of length $\leqslant
md_G(C)$.
\begin{itemize}
\item[{\rm (a)}] Is it true that $md(H)$ is
finite?

\item[{\rm (b)}] Is it true that $N^*_G(H)$ is union of finitely many double cosets
of $H$?
\end{itemize}
\end{question}

\section{Computing normal forms}

\subsection{Computing reduced forms: Algorithm I}

In this Section we discuss the standard algorithm to compute reduced
forms of elements of a group $G=A*_C B$. Suppose that the Membership
Search Problem  {\bf MSP} for $C$ is decidable in $A$ and $B$.

 Observe, that given a word $g\in F(X\cup Y)$ one
can effectively present it as a product
\begin{equation} \label{eq:dontknow}
g=g_1\cdots g_k,
\end{equation}
where $g_1, \ldots, g_k$ are reduced words in $X$ or in $Y$ and if
$g_i$ is a word in $X$ then $g_{i+1}$ is a word in $Y$ and vice
versa.

\bigskip
\noindent {\sc Algorithm I: Computing Reduced Forms.}

\medskip

\noindent {\sc Input:} a product $g=g_1\cdots g_k$ in the form
(\ref{eq:dontknow}).
 \begin{description}
\item[{\sc Step 1.}] \

\noindent
Check if $g_i \in C$, $i=1,\ldots,k$ or not. If none of the
$g_i$'s lies in $C$ then (\ref{eq:dontknow}) is reduced.
 \end{description}
\begin{description}
\item[{\sc Step 2.}] \

\noindent
Otherwise, we look at the first on the left word $g_i$
such that $g_i \in C$ and transform the word $g$ according to one of
the rules:

\begin{itemize}
 \item If $g_i \in F(X)$ then rewrite $g_i$
into a product $c_i(u_1, \ldots,u_n)$ of the given generators of
$C^\phi$, replace $g_i$ by $c_i(v_1, \ldots, v_n)$, using
substitution $u_j \rightarrow v_j, j = 1, \ldots,n$, and then
replace $g= g_1\cdots g_k$ by $g_1 \ldots g_{i-2}g_i'g_{i+2} \ldots
g_k$, where $g_i' = g_{i-1}c_i(v_1, \ldots, v_n)g_{i+1}$;

\item If $g_i \in F(Y)$ then rewrite $g_i$
into a product $c_i(v_1, \ldots,v_n)$ of the given generators of
$C^\psi$, replace $g_i$ by $c_i(u_1, \ldots, u_n)$, using
substitution $v_j \rightarrow u_j, j= 1, \ldots, n$, and then
replace $g= g_1\cdots g_k$ by $g_1 \ldots g_{i-2}g_i'g_{i+2} \ldots
g_k$, where $g_i' = g_{i-1}c_i(u_1, \ldots, u_n)g_{i+1}$;
 \end{itemize}
 thus decreasing the  syllable length
$l(g)$ of $g$. Go to Step 1.
  \end{description}

\medskip
\noindent {\sc Output:} The word
$$g = u_1\cdots u_m$$
which is  a reduced form of $g$.

\hfill {\sc End of Algorithm II}

\medskip
Notice that to carry out this algorithm one needs to be able to
verify whether or not a given  word $g_i \in F(X)$ ($g_i \in F(Y)$)
belongs to the subgroup $C^\phi$ in $A$ ($C^\psi \in B$), and, if
so, then to rewrite $g_i$ as a word in the given generators $U$ of
$C^\phi$ ($V$ of $C^\psi$). Hence, the Search Membership Problem
$\mathbf{SMP}$ has to be  decidable for the subgroup $C$ in $A$ and
$B$.

\begin{lemma}{\bf \cite{mks}}
Let $G=A*_C B$ and {\bf MSP} is decidable for $C$ in $A$ and in $B$
then Algorithm I finds a reduced form of elements of $G$.
\end{lemma}

\subsection{Computing normal forms: Algorithm II}
\label{se:can}

In this section we discuss the following algorithmic problem.

\medskip\noindent
{\bf Normal Forms Search Problem:} {\em Let $G = A \ast_C B$ and let
$S, T$ be recursive sets of representatives of $A$ and $B$
 modulo $C$. Give an  algorithm which  for a given $g \in F(X \cup Y)$
 finds the normal form of $g$ in $G$ with respect to the sets $S$ and
$T$. }

\medskip\noindent
Now  we describe the standard known decision algorithm for the
problem above (see, for example, \cite{miller1}) provided we are
given decision algorithms for the Membership Search Problem {\bf
MSP}  and the Coset Representatives Search Problem {\bf CRSP}  for
the subgroup $C$ in $A$ and in $B$ (relative to the sets $S$ and
$T$).

Given a word $g \in F(X \cup Y)$ one can effectively present it as a
product
\begin{equation}
 g = g_1\cdots g_k,
\label{eq:1}
\end{equation}
where $g_1,\ldots, g_k$ are reduced words in $X$ or in $Y$, and if
$g_i$ is a  word in $X$, then $g_{i+1}$ is a word in  $Y$, and vice
versa.

Modulo decision algorithms for problems  {\bf MSP} and  {\bf CRSP}
for the subgroup $C$ in the groups $A$ and $B$  the process of
computing the normal form is the following.

\bigskip
\noindent {\sc Algorithm II: Computing Normal Forms.}

\medskip
\noindent {\sc Input:} a word $g=g_1\cdots g_k$ in the form
(\ref{eq:1}).
 \begin{description}
\item[{\sc Step $1$.}] \

\begin{itemize}
\item [(a)] If $g_k$  is a word in $X$, then:
 \begin{itemize}
  \item  [(a.1)] Write it as
$g_k = c_ku_k$ with $u_k \in S$, $c_k \in C$ given as a word
$c_k(u_1,\ldots,u_n)$ in the given generators $U$ of $C^\phi$ (using
{\bf MSP} and {\bf CRSP}).
 \item [(a.2)] Then rewrite $c_k(u_1,\ldots,u_n)$  into
$c_k(v_1,\ldots,v_n)$.
 \item  [(a.3)] If $u_k \neq 1$ then replace $g=g_1\cdots
g_k$ by $g=g_1\cdots g_{k-2}g_{k-1}'u_k$, where $g_{k-1}' =
(g_{k-1}c_k(v_1,\ldots,v_n))$. Go to Step 2.
 \item  [(a.4)] Otherwise replace $g=g_1\cdots g_k$
by $g=g_1\cdots g_{k-2}g_{k-1}'$. Go to Step 1.
  \end{itemize}
\item [(b)] If $g_k$  is a word in $Y$ then:
  \begin{itemize}
   \item [(b.1)]  Write it as
$g_k = c_ku_k$ with $u_k \in T$ $c_k \in C$ given as a word
$c_k(v_1,\ldots,v_n)$ in the fixed generators $V$  of $C^\psi$
(using {\bf MSP} and {\bf CRSP}).
 \item [(b.2)] Then rewrite $c_k(v_1,\ldots,v_n)$
into $c_k(u_1,\ldots,u_n)$.
 \item [(b.3)] If $u_k \neq 1$ then replace
$g=g_1\cdots g_k$ by $g=g_1\cdots g_{k-2}g_{k-1}'u_k$, where
$g_{k-1}' = (g_{k-1}c_k(u_1,\ldots,u_n))$. Go to Step 2.
 \item [(b.4)] If $u_k =  1$ then replace
$g=g_1\cdots g_k$ by $g=g_1\cdots g_{k-2}g_{k-1}'$. Go to  Step 1.
 \end{itemize}

 \end{itemize}
 \end{description}

\begin{description} \item[{\sc Step 2.}] If $g$ is represented in the
form (\ref{eq:1}):
$$
g = g_1\cdots g_i u_{i+1}\cdots u_m,
$$
where $u_{i+1}\cdots u_m$ is in the  normal form, {\sc do:}
\begin{itemize}
\item [(a)] Execute Step 1 on $g_1\cdots g_i$.
  \item [(b)] If the outcome of the Step 1 on $g_1\cdots g_i$ is $g_1\cdots
  g_{i-1}'$, i.e., $u_i = 1$,  then replace $g = g_1\cdots g_i u_{i+1}\cdots u_m$
  by $g = g_1\cdots g_{i-1}' u_{i+1}\cdots u_m$ and go to Step 2.
    \item [(c)] If the outcome of the Step 1 on $g_1\cdots g_i$ is $g_1\cdots
  g_{i-1}'u_i$, i.e., $u_i \neq 1$,  then:
       \begin{itemize}
  \item [(c.1)] If both $u_i$ and $u_{i+1}$ are in $S$ (or both of them are in $T$),
rewrite $ u_iu_{i+1}$ into $c' u'_i$ where  $ c' \in C$, given as a
word $c'(v_1,\ldots,v_n)$ (or $c'(u_1,\ldots,u_n)$),   and $u'_i \in
S$ (or $u'_i \in T$).
 \item [(c.2)] If $u'_i \neq 1$ then replace $g = g_1\cdots g_i u_{i+1}\cdots
 u_m$ by $g = g_1\cdots g_{i-1}''u_i' u_{i+1}\cdots u_m$, where
 $g_{i-1}'' = g_{i-1}'c'$. Go to Step 2.
  \item [(c.3)] If $u'_i = 1$ then replace $g = g_1\cdots g_i u_{i+1}\cdots
 u_m$ by $g = g_1\cdots g_{i-1}'' u_{i+1}\cdots u_m$, where
 $g_{i-1}'' = g_{i-1}'c'$. Go to Step 2.

\end{itemize}

    \end{itemize}

 \end{description}

\noindent {\sc Output:} The word
$$g = c_1u_1\cdots u_m$$
which is  the normal form of $g$ relative to the set of
representatives $S$ and $T$.

\hfill {\sc End of Algorithm II}

\bigskip

We summarize the discussion above as the following theorem
\begin{theorem} \label{th:AlgI}
Let\/ $G = A \ast_C B$ and the problems {\bf MSP} and {\bf CRSP}
are decidable in $A$ and in $B$ for the subgroup $C$. Then Algorithm
II finds the normal form of $g$ for every given element $g \in G$.
\end{theorem}

\subsection{Complexity of Algorithm II}
\label{se:comp-Alg-I}

 Now we discuss briefly time-complexity of Algorithm II.
Recall that
 the {\em time function} $T_{\mathcal{A}}$ of an  algorithm
 $\mathcal{A}$ is defined on an input $g$ of $\mathcal{A}$ as
 the number of steps required by the algorithm $\mathcal{A}$ to
halt on the input $g$.

Obviously, the  complexity  of the time function $T_I$ of the
Algorithm II depends on complexity of the time functions of decision
algorithms for {\bf MSP} and {\bf CRSP} for  $C$ relative to $A$ and
$B$. Also, it depends on how the length of the words $c_i$ grows
during the execution of Algorithm II.

 Complexity of {\bf MSP} and {\bf CRSP}
depends on particular groups $A$, $B$, and $C$. For example, if  $A$
and $B$ are  free groups, then these problems have  linear time
complexity for a fixed subgroup $C$ (see, for example, \cite{km}).

Estimating the complexity of the rewriting process (a) is more
demanding, even in the case of amalgamated products of free groups.
Recall, that in the rewriting process (a), executing the
instructions (a.2) or (b.2),  we rewrite a word $c_{j+1}(u_1,
\ldots, u_n)$ into a word $c_{j+1}(v_1, \ldots, v_n)$. Set
$$\lambda(u,v) = \frac{\max\{|u_1|, \ldots, |u_n|\}}{\min\{|v_1|, \ldots,
|v_n|\}}$$ Then we have an upper bound
 estimate on the increase of the length
 $$|c_{j+1}(v_1, \ldots, v_n)| \leqslant \lambda(u,v) \cdot |c_{j+1}(u_1, \ldots,
 u_n)|.$$
Similarly, in the case when  we rewrite a word $c_{j+1}$ given in
the generators $v_i$ into a word in  generators $u_i$ we have an
estimate with the factor $\lambda(v,u)$. Therefore, if we denote
$$\lambda = \max\{\lambda(u,v), \lambda(v,u)\}$$
then at any rewriting step one has increase in length of at most
by the factor $\lambda$.

Now suppose, for simplicity, that the length of $c_{j}$ increases in
executing all other instructions, different from for (a.2) and
(b.2),   at most by $M + |g_j|$ where $M$ is a fixed constant (we
make this assumption to focus on the processes (a.2) and (b.2)).
Under these assumptions
 \begin{equation}
 \label{eq:step1}
 |c_j| \leqslant \lambda \cdot|c_{j+1}| + M +|g_j|
 \end{equation}
 In particular, if the length of $c_j$ does not increase at all in the rewriting
 processes other than (a.2), (b.2),  then  in $k$ steps  we will
have an exponential estimate
$$|c_1| \leqslant \lambda^{k-1}  \cdot |c_k|$$
where $k = l(g)$. So if $\lambda > 1$ then we might have
exponential growth of the length of the words $c_i$. The example
below shows that this happens in the worst case scenario.

\begin{example} Let $A = F(a,b,d), B = F(\tilde{a}, \tilde{b},
\tilde{d})$ be two free groups of ranks 3. Consider two subgroups
of rank 2:
$$C = \langle a^p, b \rangle \leqslant A,  \tilde C = \langle \tilde{a},
\tilde{b}^p \rangle \leqslant B,$$
 where $p\geqslant 2$ is an integer. Then the map
$\phi$ defined by $\phi(a^k) = \tilde a, \phi(b) = \tilde{b}^k$
gives rise to an isomorphism $\phi:C \rightarrow \tilde C$.  Put
$$G = A \ast_{C = \tilde{C} }B = \langle a,b,d,\tilde{a}, \tilde{b},
\tilde{d} \mid a^p  = \tilde a, b = \tilde{b}^p \rangle.$$
 Let $S$ be a recursive set of representatives of $A$ modulo $C$ such that the
representative in $S$ of the coset $Cda^{pm}$ is $b^{-pm}da^{pm}$
for all integers $m$. In particular,
$$da^{pm} = b^{pm}(b^{-pm}da^{pm}) \ \ \ (m \in \mathbb{Z})$$
It is not difficult to construct such $S$ since the set of
elements of the type $da^{pm}$ is recursive, as well as cosets of
$C$.  Similarly, let $T$ be a recursive set of representatives of
$B$ modulo $\tilde{C}$ such that the representative in $T$ of the
coset $\tilde{C}\tilde d \tilde{b}^{pm}$ is
$\tilde{a}^{-pm}\tilde{d} \tilde{b}^{pm}$ for all integers $m$,
which implies that
$$\tilde d \tilde{b}^{pm} = \tilde{a}^{pm}(\tilde{a}^{-pm}\tilde{d}
\tilde{b}^{pm}).$$ Now consider the following  element in $G$:
$$g = \tilde{d} d \tilde{d} d \cdots \tilde{d} d \tilde{a} = g_1 \cdots g_k$$
Then, in the notations of Algorithm II, the rewriting processes
(a.2)  and (b.2) go as follows:
$$c_k  = \tilde a = a^p$$
$$g_{k-1} = d, g_{k-1}c_k = da^p = b^p(b^{-p}da^p) = b^p u_{k-1}$$
Now the next step will be
$$c_{k-1} = b^p = \tilde{b}^{p^2}, g_{k-2} = \tilde{d}$$
Hence
$$g_{k-2}c_{k-1} = \tilde{d} \tilde{b}^{p^2}= \tilde{a}^{p^2}(\tilde{a}^{-p^2}\tilde{d}
\tilde{b}^{p^2})= \tilde{a}^{p^2}\cdot u_{k-2} = c_{k-2}\cdot
u_{k-2}.$$ In this case $\lambda = p$, lengths of the words $c_i$ do
not change in the rewriting processes other then (a.2), (b.2), so
the word $c_i$ grows every step by a factor of $p$, hence
$$|c_1| =  p^k$$ where $k = l(g)-1$.
\end{example}

\bigskip

\begin{example}
Let $A=F(a,b)$, $B=F(a',b')$ be two free groups of rank 2.
Consider two subgroups of rank 2, $C=\left<a,a^b\right>$,
$C'=\left<{a'}^p, a'^{b'}\right>$, where $p \geqslant 2$ is an
integer. Then the map $\phi$ defined by $\phi(a)= a'^{b'}$ and
$\phi(a^b)=a'^p$ gives rise to an isomorphism $\phi C \rightarrow
C'$. Put
$$
G=A*_{C=C'}B=\left<a, b, a',b'|a=a'^{b'}, a^b=a'^p \right>
$$
Let $S$ be a recursive set of representatives of $A$ modulo $B$
such that every element from $\left< b \right>$ is in $S$.
Analogously, $T$ is the set of representatives $B$ modulo $C'$
such that every element from $\left< b' \right>$ is in $T$. Now
consider the following element in $G$:
$$
g=(bb')^{-n}a(bb')^n=a^{p^n}
$$
Rewriting of this element into the normal form involves the
exponential growth of the lengths of intermediate words $c_i$.
\end{example}

Now we turn to the complexity of   rewriting processes in Algorithm
II other than (a.2), (b.2). In general, this complexity depends on
the particular algorithms for solving {\bf MSP} and {\bf CRSP} for
$C$ in $A$ and $B$. In the case of free groups $A$ and $B$ the
decision algorithm in \cite{km} for solving {\bf MSP} and {\bf CRSP}
have some important features. If we denote by $\bar w$
 the representative of the coset $Cw$ produced by the algorithm on
 the input word $w$, then the following conditions hold:

\begin{itemize}
\item For a given $w \in A$ the representative $\bar w$  of the
coset $Cw$ has the minimal possible length in $Cw$.
 \item There exists a constant $M$  such that
 for a given $w \in A$ if $w = c\bar{w}$
 for a (unique) $c \in
 C$ then $|c| \leqslant |w| + M$.
 \item the time spent by the algorithm on an input $w$  is bounded
 from above by $L|w|$ for some fixed constant $L$.
 \end{itemize}

This allows one to estimate the complexity of Algorithm II in the
case of free groups. From now on we assume that Algorithm II has
subalgorithms for solving {\bf MSP} and {\bf CRSP} which satisfy the
conditions above.

\begin{lemma}
\label{co:comp10} Let $A \ast_C B$ be a free product of free groups
with finitely generated amalgamated subgroup $C$. Then the lengths
of the words $c_i$ that occur in computations with Algorithm II on
an input $w$ is bounded from above by
 \begin{equation}
 \label{eq:length}
\lambda^k \frac{|w|+M}{\lambda -1},
\end{equation}
 where  $k =l(w)$.
\end{lemma}

\begin{proof} Let  $w = g_1 \ldots g_k$ be an input for Algorithm II in the
form (\ref{eq:1}), where $k = l(w)$. It requires $k$ steps for
Algorithm II to produce the input. According to (\ref{eq:step1}) on
each step the length of the word $c_j$   is bounded by
 $$|c_j|
\leqslant \lambda \cdot|c_{j+1}| + M +|g_j| \leqslant \lambda
\cdot|c_{j+1}| + M +|w| .$$
  Hence in $k$ steps we will have the following estimate on
  the lengths of the words $c_j$, $j = 1, \ldots k$.
\begin{eqnarray*} \lambda (\cdots (\lambda(|w|+M) + |w|+M))\cdots ) + |w| +M
&=& (\lambda^{k-1} + \cdots +1)(|w|+M)\\ &\leqslant&  \lambda^k
\frac{|w|+M}{\lambda -1}, \end{eqnarray*}
 as required.
\end{proof}

\begin{corollary}
\label{co:comp10b} Let $A \ast_C B$ be a free product of free groups
with finitely generated amalgamated subgroup $C$. Then the time
spent by Algorithm II on an input $w$ is bounded above by
$$k  \cdot L_1 \cdot |w| \cdot \lambda^k \cdot(|w|+M)$$
where $L_1$ is a fixed constant and $k = l(w)$.
\end{corollary}

\begin{proof} Indeed, Algorithm II works $k$ steps on an input $w$ with $l(w)
= k$. On each step it rewrites a current word $c_j$ of the length
bounded from above in (\ref{eq:length}). The rewriting involves the
subalgorithms  for solving {\bf MSP} and {\bf CRSP}. These
algorithms spend at most linear time with respect to the length of
the input. Putting all the estimates together we have the
 resulting estimate above.
\end{proof}

 Combining the corollaries above with the example we
have the following result.

\begin{theorem}
\label{th:comp} \
\begin{enumerate}
\item [{\rm (1)}] Let $A \ast_C B$ be a free product of free groups with
finitely generated amalgamated subgroup $C$. Then Algorithm II has
at most exponential (in the length of the input words) time
complexity  function bounded by:
$$k L_1 |w| \lambda^k (|w|+M)$$
where $k, L_1, \lambda, M$, and $w$ are as above;
 \item [{\rm (2)}] There are finitely generated free groups  $A$ and
$B$ and a finitely generated  subgroup $C$ in $A$ and $B$ such that
in the free product with amalgamation $A \ast_C B$ the Algorithm II
has precisely the   exponential time complexity as above.
\end{enumerate}
\end{theorem}

However, we will show  in the subsequent paper \cite{amalgam2} that
the situation in the example above is very rare, and in every free
product with amalgamation $G = A \ast_C B$ of free groups with a
finitely generated group $C$ the Algorithm II is  very fast on
generic inputs.

\subsection{Computing cyclically reduced normal forms: Algorithm III}
\label{subsec:3.3}

  In this section, we shall discuss the standard algorithm to
find a  cyclically reduced normal form of an element $g$ of a group
$G = A \ast_C B$. As before, we assume that the element $g$ is
given in the form (\ref{eq:1}):
$$g = g_1\cdots g_n,$$ where $g_1,\ldots, g_n \in F(X) \cup F(Y)$
 and  $g_i \in F(X)$ if and only if $g_{i+1} \in F(Y)$.

  We work under assumption that the Membership Search Problem {\bf MSP},
  the Coset Representative Search Problem
{\bf CRSP},  and the Conjugacy Membership Search Problem {\bf CMSP}
are decidable in $A$ and $B$ for the subgroup $C$, and
 we have the decision algorithms in our
possession. Notice that we need {\bf CMSP} only because we have a
slightly stronger notion of reduced forms than the usual one (see
Section \ref{sec:1.2}).

 Observe, that the uniform version of {\bf CMSP} is decidable in free groups
 and the  decision
 algorithm has linear time complexity (in the length of the input word $w$) for a given
 finitely generated subgroup $C$  \cite{km}.

\bigskip
\noindent {\sc Algorithm III: Computing Cyclically Reduced Forms.}

\medskip
\noindent {\sc Input:} a word $g$ in the form (\ref{eq:1}).
 \begin{description}
\item [{\sc Step $1$}] \ Find the normal  form of $g$ using the
Algorithm II:
$$g = c_gg_1\cdots g_k.$$
Observe that $l(g) = k$ and for every $g_i$ we know its factor
$F(g_i)$.

 \item [{\sc Step 2}] \

 \begin{itemize}
\item [(a)] If $l(g)=0$ then $g$ is already in cyclically reduced from.

\item [(b)] If $l(g)=1$, for example, if $g \in A$, then  check whether
$g$ is a conjugate of  an element $c\in C$ or not, using the
algorithm for {\bf CMSP}.    In the former case, $c$ is a
cyclically reduced form of $g$ and  the algorithm for {\bf CMSP}
gives one of such elements $c$. In the latter case, $g$ is already
in cyclically reduced form.

\item  [(c)] Let $l(g) \geqslant 2$.

\begin{itemize} \item If $F(g_1) \neq F(g_k)$, then $g$ is already in a
cyclically reduced normal form.

\item If  $F(g_1) = F(g_k)$.  Then $g$ is conjugate to
$$ (g_kc_qg_1)g_2\cdots g_{k-1}.$$
 Now apply the decision algorithms for {\bf MSP} and {\bf CRSP} to the word $(g_kc_qg_1)$
 to find the normal form $c'g'_1$ of it.  If $g'_1 \neq 1$ then
$$c'g'_1g_2\cdots g_{k-1}$$
is a cyclically reduced normal form of $g$. Otherwise,
$$F(c'g_2) = F(g_{k-1})$$ and we apply  the procedure above to
$c'g_2\cdots g_{k-1}.$ \end{itemize} \end{itemize} \end{description}
\hfill {\sc End of Algorithm III}

\begin{lemma}
\label{le:modified-III}
 Let\/ $G = A \ast_C B$ and the problems {\bf
MSP} and {\bf CRSP} are decidable in $A$ and $B$ for the subgroup
$C$. Then there exists an algorithm that for a given element $g \in
G$ finds an element $g' \in G$ such that $g'$ is a conjugate of $g$
and if $l(g')
> 1$ then  $g'$ is a cyclically reduced normal
  form of $g$.
\end{lemma}
  \begin{proof}
   Direct analysis of Algorithm III shows that a decision algorithm for the  problem {\bf
   CMSP} is used only when executing instructions in the case (b).
   However, if we modify Algorithm III in such a manner that it stops immediately when
   the case (b) occurs, then the modified algorithm satisfies the requirements of the lemma.
     \end{proof}

\begin{theorem} \label{th:AlgII}
Let\/ $G = A \ast_C B$ and the problems {\bf MSP}, {\bf CRSP}, and
{\bf CMSP} are decidable in $A$ and $B$ for the subgroup $C$.   Then
for a given element $g \in G$  Algorithm III finds a cyclically
reduced normal form of $g$ in  time $T_{III}(g)$ which can be
bounded from above as follows:
 $$T_{III}(g) \leqslant T_{II}(g) + K \cdot \max\{T_{CMSP}(c_gg_1),
 \left (T_{MSP}(g_kc_gg_1) + T_{CRSP}(g_kc_gg_1)\right )\cdot l(g)\},$$
where $T_{II}$, $T_{SMP}$, $T_{CMSP}$, $T_{CRSP}$ are the time
functions, correspondingly, of  Algorithm II, and the decision
algorithms for ${\bf MSP, CMSP, CRSP}$, $K$ is a constant.  In
particular, if $A$ and $B$ are free groups then
$$T_{III}(g) \leqslant T_{II}(g) + K_1 \cdot |g|\cdot l(g),$$
where $K_1$ is a constant\/ {\rm (}depending on $C${\rm )}  and\/
$|g|$ is the length of the input $g$ given as a word in $F(X \cup
Y)$.
\end{theorem}

\bigskip

\section{Regular Elements and black holes} \label{se:5-1}

\subsection{Bad pairs}

Let $G = A \ast_C B$.

\begin{definition}
  We say that
$(c,g) \in C \times G$ is a bad pair if $ c \neq 1, g \not \in C$,
and $gcg^{-1} \in C$.
\end{definition}

Notice that if $(c, g)$ is a bad pair then $g \in N^\ast_G(C)
\smallsetminus C$ and $c \in Z_g(C)$. The following lemma gives a
more detailed description of bad pairs.

\begin{lemma}
\label{le:B}
 Let $c \in C \smallsetminus \{1\}$ and $g \in G \smallsetminus C$.
 If $g = c_gp_1 \cdots p_k$ is the
 normal  form of $g$ then
$(c,g)$ is a bad pair if and only if the following system
$B_{c,g}$ has a solution $c_1, \ldots, c_{k}$ with $ c_i \in C$:
\begin{eqnarray*}
p_kc p_k^{-1} & =  & c_{1}\\
p_{k-1} c_{1}p_{k-1}^{-1} & =& c_{2} \\
 & \vdots &  \\
p_1c_{k-1} p_1^{-1} & =& c_k \label{ar:bp}
\end{eqnarray*} Moreover, in this case $p_i \in N^\ast_{F(p_i)}(C)$ and
$c \in Z_A(C) \cup Z_B(C)$. \end{lemma}

\begin{proof}  This lemma is a particular case of Lemma \ref{le:systems}.
\hfill \end{proof}

Observe, that consistency  of the system $B_{c,g}$ does not depend
on a particular choice of  representatives of $A$ and $B$ modulo
$C$. Sometimes we shall treat $c$ as a variable, in which case the
system will be denoted $B_g$. If $c,c_1, \ldots, c_{k+1} \in C
\smallsetminus \{1\} $ is a solution of $B_g$  then we call it a
{\em nontrivial} solution of $B_g$.

 Now we will study slightly more general equations of the type $gc
 = c'g'$ and their solutions $c, c'$ in $C$.

\begin{lemma}
\label{le:systems} Let $G = A \ast_C B$,  $g, g^\prime \in G$ be
elements given by their  normal forms:
\begin{equation}
g = c_{g} p_1 \cdots p_k, \ \ \ g' = c_{g'} p'_1 \cdots p'_k  \ \
\ (k \geqslant 1). \label{eq:4}
\end{equation}
 Then the equation $gc = c'g'$ has a solution $c, c' \in C$ if and only
if the following system $S_{g,g'}$ in variables $c, c',  c_1,
\ldots, c_{k}$ has a solution in $C$:
\begin{eqnarray*}
p_kc & =  & c_1p'_k\\
p_{k-1} c_1 & = & c_2 p'_{k-1}\\
 & \vdots &  \\
p_1c_{k-1} & =& c_{k}p'_1\\
c_g c_k& = & c' c_{g'}
\end{eqnarray*}
\end{lemma}
{\it Proof.} Let  $c, c' \in C$ be  a solution to the equation
$gc= c'g'$. We then rewrite the  equality $gc= c'g'$ as
$$
c_{g} p_1 \cdots p_k c= c' c_{g'}p'_1 \cdots p'_k.
$$
Notice that the right hand side of this equality is in the normal
form. Following Algorithm II we shall rewrite the left hand side of
this equality  into the normal form. After rewriting the both sides
must coincide as the normal  forms of the same element. This gives
rise to the system of equations for some elements  $c, c', c_1,
\ldots, c_{k} \in C$, as above. Conversely, if the system $S_{g,g'}$
has a solution then the elements $c, c'$ give a solution of the
equation $gc= c'g'$.
 \hfill$\Box$

The first $k$ equations of the system $S_{g,g'}$  form what we
call the {\em principal system} of equations, we denote  it by
$PS_{g,g'}$. In what follows  we  consider  $PS_{g,g'}$ as a
system in variables $c, c_1, \ldots, c_k$ which take values in
$C$,  the elements $p_1, p'_1, \ldots, p_k, p'_k$ are constants.

\subsection{Regular elements}

 Now we specify, in our particular context, the general concepts of
a ``black hole" and ``regular part'' as discussed in the
Introduction.

\begin{definition}
The set
$$BH = ({N^\ast}_G(C) \smallsetminus C) \cup  {Z}_G(C) $$
 is called a {\em black hole}. Elements from $BH$ are called
 {\em singular}, and elements from $G \smallsetminus BH$ {\em regular}.
\end{definition}
Notice that if  the subgroup $C$ has a finite malnormality degree in
$G$  then every element $g$ with $l_0(g) > md_G(C)$  is regular. In
particular, it follows from Lemma \ref{le:mal-deg} that if $C$ is
malnormal in $A$ or in $B$ then every  element $g \in G$ with $l(g)
\geqslant 2$ is regular. Notice also, that if $g \in G
\smallsetminus C$ is regular then all elements in $CgC$ are regular.

Observe, that the condition 1) in the Conjugacy Criterion, indeed,
does not apply for regular  elements.

 The following description of singular elements
follows from Lemma \ref{le:B}.
\begin{corollary}
\label{co:DescReg}
 Let $G = A \ast_C B$. Then:
\begin{enumerate}
 \item [1)]  an element $g \in G \smallsetminus C$  is singular if
and only if the system $B_{g,c}$  has a solution $c, c_1, \ldots,
c_{k}$, where $c, c_i$ are non-trivial elements from $C$;
 \item [2)] $Z_G(C) = Z_A(C) \cup Z_B(C)$
\end{enumerate}
\end{corollary}

As we have seen already, an element $g \in G$ is singular if
 and only if the system $gc = c_1g$ has a nontrivial solution $c,
 c_1$ in $C$.

\subsection{Effective recognition of regular elements}

\begin{definition}
\label{de:shifts} Let $M$ be a subset  of a group $G$. If\/ $u,v
\in G$ then the set $uMv$ is called a $G$-shift of $M$. For a set
${\mathcal{M}}$ of subgroups of $G$ denote by
$\Phi({\mathcal{M}},G)$ the least set of subsets of $G$ which
contains ${\mathcal{M}}$ and is closed under $G$-shifts and
intersections.
\end{definition}

\begin{lemma}
\label{le:5b} Let $G$ be a group and $C$ be a subgroup of $G$. If
$D \in \Phi(\{C\},G)$, $D \neq \emptyset$ then there exist
elements $ g_1, \dots, g_n, h \in G$ such that
$$
D=(C^{g_1}\cap \cdots \cap C^{g_n})h
$$
In particular, non-empty sets in $\Phi(\{C\},G)$ are particular
cosets from $G$.
\end{lemma}

\begin{proof} Induction on the number of operations required to construct
$D$ from $C$.
 For a tuple $\bar g = (g_1, \ldots,
 g_n)$ of elements  from $G$ put
$$C_{\bar{g}}=(C^{g_1}\cap \dots \cap C^{g_n}).$$
  Let $D = C_{\bar{g}}h$ for some $\bar{g} \in G^n, h \in G$. Then
  for any $a, b \in G$:
 $$a D b =D^{a^{-1}}a b = C_{\bar{g}a^{-1}} ab,$$
where $\bar{g}a^{-1} = (g_1a^{-1}, \ldots, g_na^{-1})$, i.e.,
$aDb$ is in the required form.

Observe, that for arbitrary subgroups $K, L \leqslant G$ and
elements $a, b \in G$ if $h \in Ka \cap Lb$ then
 \begin{equation}
 \label{eq:Cosets}
 Ka \cap Lb = (K \cap L)h.
 \end{equation}
Therefore, if
 $h_3 \in C_{\bar{g_1}}h_1
\cap C_{\bar{g_2}}h_2$, then
$$C_{\bar{g_1}}h_1
\cap C_{\bar{g_2}}h_2 = (C_{\bar{g_1}} \cap C_{\bar{g_2}})h_3 =
C_{\bar{g_3}}h_3,$$ where  $\bar{g_3}$ is concatenation of
$\bar{g_1}$ and $\bar{g_2}$. \end{proof}

\begin{lemma}
\label{le:PS} Let $G = A \ast_C B$. Then for given two elements $g,
g' \in G$ in their normal forms
 $$g = c_{g} p_1 \cdots p_k, \\ \ g' = c_{g'} p'_1 \cdots p'_k   \ \ \ (k \geqslant 1)$$
 the set $E_{g,g'}$, of all elements $c$ in $C$ for which the system $PS(g,g')$ has a solution
 $c, c_1, \ldots, c_k \in C$, is equal to
$$E_{g,g'} = C \ \cap p_k^{-1}Cp'_k \ \cap \ \cdots \ \cap p_k^{-1}\cdots p_1^{-1} C p'_1\cdots p'_k.$$
In particular, if $E_{g,g'} \neq \emptyset$ then $E_{g,g'} =
C_{g,g'}c_{g,g'}$ for some subgroup $C_{g,g'} \leqslant C$ and
some element $c_{g,g'} \in C$.
\end{lemma}

\begin{proof} Let
$$g = c_{g} p_1 \cdots p_k, \ \ \ g' = c_{g'} p'_1 \cdots p'_k.$$
Denote by $V_i$ the set of all solutions $(c, c_1, \ldots,c_i) \in
C^{i+1}$ of the system formed by the first $i$ equations of
$PS(g,g')$. Let $D_{m,i}$ be the projection of $V_i$ onto its
$m$-s component.

The first equation of the system $PS(g,g')$ gives:
$$p_kc_0 (p'_k)^{-1}= c_1$$
where for uniformity we denote $c$ by $c_0$. Therefore,
$$D_{1,1} = p_kC(p'_k)^{-1} \cap C, \ \  D_{0,1} = p_k^{-1}D_{1,1}p'_k$$
and  $(c_0,c_1) \in V_1$  if and only if
$$c_1 \in D_{1,1}, \ \ c_0 = p_k^{-1} c_1 p'_k.$$
Clearly, the sets $D_{0,1}$ and $D_{1,1}$   are in $\Phi_C$.

Now we rewrite the $i$-s equation $p_{k-i+1}c_{i-1}  =
c_{i}p'_{k-i+1}$ of the system $PS(g,g')$ in the form
$$p_{k-i+1}c_{i-1}(p'_{k-i+1})^{-1}  = c_{i}$$
It follows that
 \begin{equation}
 \label{eq:Dii)}
 D_{i,i} = p_{k-i+1}D_{i-1,i-1}(p'_{k-i+1})^{-1}
\cap C,
\end{equation}
 where $i = 1, \ldots ,k$ and $ D_{0,0} = C$. In
particular
$$D_{k,k} = p_1D_{k-1,k-1}(p'_1)^{-1} \cap C$$
Clearly, $(c, c_1, \ldots, c_k)$ is a solution of the system
$PS(g,g')$ if and only if $c_k \in D_{k,k}$ and  $c_{i-1} = p_{k-i
+1}^{-1}c_ip'_{k-i+1}$.  More precisely, since

  $$D_{i-1,k} = p_{k-i +1}^{-1}D_{i,k}p'_{k-i+1}$$
 it follows now that,


$$D_{k-i,k} = D_{k-i,k-i} \ \cap p_i^{-1}Cp'_i \ \cap \ \cdots \
\cap p_i^{-1}\cdots p_1^{-1} C p'_1\cdots p'_i.$$
 In particular,
 $$E_{g,g'} = D_{0,k} = C \ \cap p_k^{-1}Cp'_k \ \cap \ \cdots \
\cap p_k^{-1}\cdots p_1^{-1} C p'_1\cdots p'_k.$$ So $E_{g,g'} \in
\Phi(C,G)$. By Lemma \ref{le:5b}
 $$p_k^{-1}Cp'_k \ \cap \ \cdots \ \cap p_k^{-1}\cdots p_1^{-1} C p'_1\cdots
 p'_k = Hu$$
 for some subgroup $H \leqslant G$ and $u \in G$. Now we can see from
 (\ref{eq:Cosets}) that
$$E_{g,g'} = C \cap Hu = C_{g,g'}c_{g,g'}$$
 for some subgroup $C_{g,g'} \leqslant C$ and  $c_{g,g'} \in C$, as required.
\end{proof}

Denote buy $Sub(C)$ the set of all subgroups of $C$. By Lemma
\ref{le:5b} non-empty sets from $\Phi(Sub(C),A)$ (respectively,
from $\Phi(Sub(C),B)$) are some cosets of subgroups from $A$
(respectively, from $B$).

\begin{corollary}
\label{co:4} Let $G = A \ast_C B$.   If the Cardinality Search
Problem  is decidable for $\Phi(Sub(C),A)$ in $A$ and for
$\Phi(Sub(C),B)$ in $B$
 then given $g, g'$ as above, one can effectively find the
set $E_{g,g'}$. In particular, one can effectively  check whether
or not $E_{g,g'}$ is empty, singleton, or infinite.
\end{corollary}

\begin{proof} In notations of Lemma \ref{le:PS}
$$E_{g,g'} = p_k^{-1}\cdots p_1^{-1} D_{k,k} p'_1\cdots p'_k.$$
Therefore it suffices to solve the cardinality problem for the set
$D_{k,k}$. The quality \ref{eq:Dii)}
 $$ D_{i,i} = p_{k-i+1}D_{i-1,i-1}(p'_{k-i+1})^{-1} \cap C,$$
and Lemma \ref{le:5b}  show that each $D_{i-1,i-1}$ is a coset of
the type $C_ic_i$ where $C_i \leqslant C$ and $c_i \in C$. Moreover,
since the Cardinality Search Problem  is decidable for
$\Phi(Sub(C),A)$ in $A$, and for $\Phi(Sub(C),B)$ in $B$, the
equality (\ref{eq:Cosets}) shows how one can effectively find the
element $c_i$ and the direct expression for the subgroup $C_i$ (in
terms of shifts and intersections). Therefore, in $k$ steps one can
find $D_{k,k}$, and hence the set $E_{g,g'}$. Moreover, on each step
one can find the cardinality of the set $D_{i,i}$. This proves the
corollary. \end{proof}

\begin{lemma}
\label{le:PS2} Let $G = A \ast_C B$ and  $g, g' \in G $. If $l(g)
= l(g') \geqslant 1$ and  the system $PS(g,g')$ has more then one
solution in $C$ then the elements $g, g'$ are singular.
\end{lemma}

\begin{proof} Let  $c, c_1, \ldots, c_{k}$  and $b, b_1, \ldots, b_{k}$ be
two distinct  solutions of the principal system $PS(g,g')$. Denote
for uniformity $c_0 = c, b_0 = b.$ Hence we have the following
systems of equations:
$$
p_kc_{0}  =   c_{1}p'_k,  \ \ \ \  p_k b_{0}  =  b_{1}p'_k $$

$$p_{k-1} c_{1}  = c_{2} p'_{k-1},  \ \ \ \ p_{k-1} b_{1}  = b_{2}
p'_{k-1}$$
 $$  \vdots $$
  $$p_1c_{k-1} = c_{k}p'_1,  \ \ \ \
p_1 b_{k-1}  = b_{k}p'_1$$ Expressing $p'_k$ from the first two
equations in the system above, and then $p'_{k-1}$  from the next
two equations,  and so on, we get the following equalities:
 $${c_1}^{-1}p_kc_0 = {b_1}^{-1}p_kb_0$$
$${c_2}^{-1}p_{k-1}c_1 = {b_2}^{-1}p_{k-1}b_1$$
 $$ \vdots $$
$${c_k}^{-1}p_1c_{k-1} = {b_k}^{-1}p_1b_{k-1}$$
Rewriting these equalities we obtain:
$$
p_k^{-1}b_1c_{1}^{-1}p_k = b_0c_0^{-1},
$$
$$
p_{k-1}^{-1}b_2c_{2}^{-1}p_{k-1} = b_1c_1^{-1},
$$
$$ \vdots $$
$$
p_{1}^{-1}b_kc_{k}^{-1}p_{1} = b_{k-1}c_{k-1}^{-1}.
$$
Observe that all the elements $b_ic_i^{-1}$ are non-trivial. By
Lemma \ref{le:B} the element $g$ is singular. Similar argument shows
that $g'$ is also singular. \end{proof}

\medskip

The next result shows that one can effectively determine whether a
given element $g \in G $ is regular or not.

\begin{theorem}
\label{le:BHeff} Let $G = A \ast_C B$ be a free product of finitely
presented groups $A$ and $B$ amalgamated over a finitely generated
subgroup $C$. Assume that the following algorithmic problems are
decidable:

\begin{itemize}
 \item The Search Membership Problem for the subgroup $C$ in $A$ and in $B$.
  \item The Coset Representative Search Problem for the subgroup $C$ in $A$ and in $B$.
   \item The Cardinality Search Problem for $\Phi(Sub(C),A)$ in $A$ and for
   $\Phi(Sub(C),B)$ in $B$.
    \item The Membership Problem for $N^\ast_A(C)$ and  $Z_A(C)$ in $A$, and for
    $N^\ast_B(C)$ and     $Z_B(C)$ in $B$.
\end{itemize} Then there exists an algorithm to determine whether a given
element in $G$ is regular or not.
\end{theorem}

\begin{proof} For a given $g \in G$ we can find the  normal  form
of $g$ using Algorithm II. Now there are two cases to consider.

1) If $l(g) > 1$ then by Lemma \ref{le:B} $g$ is a singular
element if  and only if the system $B_{c,g}$ has a nontrivial
solution $c, c_1, \ldots, c_k \in C$. Observe, that if the system
$B_{c,g}$ has two distinct solutions then one of them is
non-trivial (i.e., $c, c_1 \ldots, c_k \neq 1$).

Now if $B_{c,g}$ has no solutions in $C$ (and we can check it
effectively) then $g$ is regular. If $B_{c,g}$ has precisely one
solution then we can find it and check whether it is trivial or
not, hence we can find out whether $g$ is regular or not. If
$B_{c,g}$ has more then one solution (and we can verify this
effectively) then $g$ is not regular.

2) If $l(g) = 1$ then $g \in A \cup B \smallsetminus C$. In this
case $g$ is regular if and only if $g \not \in N^\ast_B(C) \cup
N^\ast_B(C)$.  Since the sets $N^\ast_A(C)$ and $N^\ast_B(C)$ are
recursive one can algorithmically check  if $g$ is regular or not.

3) If $l(g) = 0$ then $g$ is regular if and only if $g \not \in
Z_G(C)$. By Corollary \ref{co:DescReg} $Z_G(C) = Z_A(C) \cup
Z_B(C)$. Since the  sets $Z_A(C)$ and $Z_B(C)$ are recursive one can
check whether or not $g$ is regular. This proves the theorem.
\end{proof}

\begin{corollary}
\label{co:free-reg-rec}
 Let $G = A \ast_C B$ be a free product
with amalgamation of free groups $A, B$. Then the set of regular
elements in $G$ is recursive.
 \end{corollary}

\begin{remark}
{\rm The decision algorithm for checking whether a given element is
regular or not is fast  ``modulo'' Algorithm II and the algorithm
${\mathcal{B}}$  for finding cardinality of sets of the type
$E_{g,g'}$. In general, both Algorithm II and ${\mathcal{B}}$ can be
 exponential in the worst case. However, we will show later that
generically both the algorithms are fast. }\end{remark}

One can improve on Theorem \ref{le:BHeff} in the following way.
Denote by $CR$ the set of all elements in $G$ which have at least
one regular cyclically reduced normal form of length greater than 1,
i.e., $CR$ is the set of elements in $G$ which are conjugates of
cyclically reduced regular elements. Now by $CR_{>1}$ we denote a
subset of $CR$ consisting  of elements of cyclically reduced length
grater than 1, so $CR_{>1}$ is the set of elements in $G$ which are
conjugates of cyclically reduced regular elements of length greater
than 1.

\begin{corollary}
\label{le:CR} Let $G = A \ast_C B$. Assume that the following
algorithmic problems are decidable:

\begin{itemize}
  \item The Search Membership Problem for the subgroup $C$ in $A$ and  $B$.
  \item The Coset Representative Search Problem for the subgroup $C$ in $A$ and  $B$.
   \item The Cardinality Search Problem for $\Phi(Sub(C),A)$ in $A$ and for
   $\Phi(Sub(C),B)$ in $B$.
    \item The Membership Problem for $N^\ast_A(C)$ and $Z_A(C)$ in $A$ and for
    $N^\ast_B(C)$ and    $Z_B(C)$ in $B$.
\end{itemize} Then  there exists an algorithm that for a given element $g \in G$ decides
whether $g$ belongs to  $CR_{>1}$  or not, and if so, then finds a
regular cyclically reduced normal form of $g$.
\end{corollary}

\begin{proof} Let $g \in G$. By Lemma \ref{le:modified-III} one
can effectively find an element $g' \in G$ such that $g'$ is a
conjugate of $g$ and if $l(g') >1$ then $g'$ is a cyclically reduced
normal form of $g$. It follows that if $l(g') \leq 1$ then $g \not
\in CR$. Suppose $l(g') >1$.  We claim that in this case   $g$ has a
cyclically reduced regular normal form, say $g_1$, if and only if at
least one of the cyclic permutations of $g'$ is regular. Indeed,
observe that $g'$ and $g_1$ are conjugated in $G$, hence by the
conjugacy criterion $g_1 = \pi_i(g')^c$ for some $i$-cyclic
permutation $\pi_i(g')$ of $g'$ and some $c \in C$. Since $g_1$ is
regular this implies that $\pi_i(g') =  g_1^{c^{-1}}$ is also
regular (easy calculation). It follows that one of cyclic
permutations of $g'$ is regular. Now one can effectively list all
cyclic permutations $\pi_j(g')$ of $g'$ and apply the decision
algorithm from Theorem \ref{le:BHeff} to each  cyclic permutation
$\pi_j(g')$ to verify if there is a regular one among them. This
proves the result.

\end{proof}

Denote by $CR_0$ a  subset of $CR$ consisting  of elements of
cyclically reduced length 0.

\begin{lemma}
\label{le:reg-C} Let $G = A \ast_C B$ be a free product of finitely
presented groups $A$ and $B$ amalgamated over a finitely generated
subgroup $C$.  Assume that the following algorithmic problems are
decidable:

\begin{itemize}
  \item The Search Membership Problem for $C$ in $A$ and  $B$.
  \item The Coset Representative Search Problem for $C$ in $A$ and  $B$.
  \item The Conjugacy Membership Search Problem for $C$ in $A$ and
  $B$.
   \item The Membership Problem for $Z_A(C)$ and $Z_B(C)$.
  \end{itemize}
  Then  there exists an algorithm that for a given element $g \in G$ decides
whether $g$ belongs to  $CR_0$  or not, and if so, then finds a
regular cyclically reduced normal form of $g$.
\end{lemma}
 \begin{proof}
Let $g \in G$.  By Theorem \ref{th:AlgII} one can use Algorithm III
to find a cyclically reduced normal form of $g$. So we may assume
from the beginning that $g$ is already in a cyclically reduced
normal form. If $l(g)
>0$ then $g \not \in CR_0$. Suppose $l(g) = 0$, i.e., $g \in C$.
  We claim that $g \in CR_0$ if and only if it is  regular.
Indeed, by conditions of the lemma $g$ is a conjugate of some
regular cyclically reduced element $g'$ which must be in $C$ (since
it has length  0 as $g$). But then since $g'$ is regular the
elements $g$ and $g'$ are conjugates in $C$. Hence  $g$ is regular
since  $g'$ is regular. Observe, that $g \in C$ is regular if and
only if $ g \not \in Z_G(C)$.  By Corollary \ref{co:DescReg} $Z_G(C)
= Z_A(C) \cup Z_B(C)$.  Now, since the Membership problem for for
$Z_A(C)$ and $Z_B(C)$ is decidable in $A$ and $B$, one can check
algorithmically if $g$ belongs to $Z_G(C)$ or not, thus solving the
question if $g$ is in $CR_0$ or not.

 \end{proof}

\begin{remark} Corollary \ref{le:CR} claims that the Search
Membership Problem for the set $CR_{>1}$ is decidable in $G = A
\ast_C B$ under corresponding conditions on the factors $A$, $B$.

\end{remark}

\section{The Conjugacy Search Problem and regular elements } \label{2.3}
\label{se:4} \label{se:5-2}

In this section we study the Conjugacy Search Problem in a group $G
= A \ast_C B$.

We start with the following particular case of the the Conjugacy
Search Problem.

\begin{theorem}
\label{th:main}
 Let $G = A \ast_C B$ be a free product of finitely
presented groups $A$ and $B$ amalgamated over a finitely generated
subgroup $C$.  Assume that the following algorithmic problems are
decidable:
\begin{itemize}
  \item The Word Problem in $A$ and in $B$.
  \item The Search Membership Problem for the subgroup $C$ in $A$ and in $B$.
  \item The Coset Representative Search Problem for the subgroup $C$ in $A$ and in $B$.
  \item The Cardinality Search Problem for $\Phi(Sub(C),A)$ in $A$ and for
   $\Phi(Sub(C),B)$ in $B$.
\end{itemize}

Then the Conjugacy Search Problem in $G$ is decidable for all pairs
from $CR_{>1} \times G$.
\end{theorem}

\begin{proof} Let $g \in CR_{>1}$ and $h \in G$. By Corollary \ref{le:CR} one can
find a regular cyclically reduced normal form $g'$ of $g$.
Meanwhile, by Lemma \ref{le:modified-III} one can find and element
$h' \in G$ such that $h'$ is a conjugate of $h$ and if $l(h')
> 1$ then  $h'$ is a cyclically reduced normal  form of $h$. Since
$g \in CR_{>1}$ its cyclically reduced length is greater then $1$,
hence if $l(h') \leq 1$ then $h$ is not a conjugate of $g$. Suppose
now that $l(h') > 1$, in this case $h'$ is a cyclically reduced
normal form of $h$. This shows that we may assume  from the
beginning that $g$ is regular   and $g$, $h$ are given in cyclically
reduced normal forms:

 $$ g = cp_1 \ldots p_k, \ \ \   h = c^{\prime}p_1^{\prime}
 \ldots  {p_{k^{\prime}}}^{\prime}.$$
 According to the conjugacy criterion, the elements $g$
and $h$ are conjugate in $G$ if and only if $k = k^\prime$ and for
some cyclic permutation $\pi(h)$ of  $h$ the equation $c^{-1}gc =
\pi(h)$ has a solution $c$ in $C$.  By Lemma \ref{le:systems} the
equation $c^{-1}gc = \pi(h)$ has a solution in $C$ if and only if
the system $S_{g, \pi(h)}$ has a solution in $C$.  Since $g$ is
regular the system $PS_{g, \pi(h)}$ has at most one solution in $C$.
Decidability of the Cardinality Search Problem problems for
$\Phi(Sub(C),A)$ in $A$ and for   $\Phi(Sub(C),B)$ in $B$  allows
one to check whether $PS_{g, \pi(h)}$ has a solution in $C$ or not,
and if it does, one can find the solution. Now one can verify
whether this solution satisfies the last equation of the system
$S_{g, \pi(h)}$ or not (using decidability of the word problem in
$A$ and $B$). If not, the system $S_{g, \pi(h)}$ has no solutions in
$C$, as well as the equation $c^{-1}gc = \pi(h)$. Otherwise, the
system $S_{g, \pi(h)}$  and the equation $c^{-1}gc = \pi(h)$ have
solutions in $C$ and we have found one of these solutions. This
proves the theorem.
\end{proof}

Now we study conjugacy search problem for regular elements of length
$\leqslant 1$.

\begin{lemma}
\label{le:reg-0} Let $G = A \ast_C B$ be a free product of finitely
presented groups $A$ and $B$ amalgamated over a finitely generated
subgroup $C$.  Assume that the following algorithmic problems are
decidable:

\begin{itemize}
  \item The Search Membership Problem for the subgroup $C$ in $A$ and in $B$.
  \item The Coset Representative Search Problem for the subgroup $C$  in $A$ and in $B$.
  \item The Conjugacy Membership Search Problem for $C$ in $A$ and
  $B$.
  \item The Conjugacy Search Problem in $C$.
  \end{itemize}
   Then the Conjugacy Search Problem in $G$ is decidable for all pairs
from $CR_0 \times G$.
\end{lemma}

\begin{proof} Let $(g,h) \in G \times G$. Since the Membership Search
Problem for the set $CR_0$ is decidable (Theorem \ref{le:reg-C}) one
can algorithmically check if $g \in CR_0$, and if so, find a regular
element $g'$ given in a cyclically reduced normal form which is a
conjugate of $g$.   By Theorem \ref{th:AlgII} one can apply
 Algorithm III to find a cyclically reduced normal form  of
$h$. Clearly, $h$ is a conjugate of $g$ if and only if $h'$ is a
conjugate of $g'$. Replacing $(g,h)$ by $(g',h')$ if necessary we
may assume that $g$ and $h$ are already in their cyclically reduced
normal forms and $g$ is regular.  It follows that $g \in C$ since $g
\in CR_0$. By the conjugacy criterion, if $h \not \in C$ then $h$ is
not a conjugate of $g$. If $h \in C$ then by the conjugacy criterion
there is a sequence of elements\/ $g =c_0,c_1,\dots, c_t=h$, where
$c_i \in C$ and adjacent elements $c_i$ and $c_{i+1}$, $i = 0,\dots,
t-1$, are conjugate in $A$ or in $B$. Since $g = c_0$ is regular -
it does not belong to $Z_G(C)$. By Corollary \ref{co:DescReg}
$Z_G(C) = Z_A(C) \cup Z_B(C)$. This implies that the element that
conjugates $c_0$ into $c_1$ must be in $C$. Hence, $c_1$ is also
regular. By induction on $t$, all the pairs  $c_i, c_{i+1}$ are
conjugated in $C$, as well as elements $g$ and $h$. Thus the Search
Conjugacy Problem for $g$ and $h$ in $G$ is reduced to the Search
Conjugacy Problem for $g$ and $h$ in $C$ - which is decidable.

 \end{proof}

\begin{lemma}
\label{le:reg-1} Let $G = A \ast_C B$ be a free product of finitely
presented groups $A$ and $B$ amalgamated over a finitely generated
subgroup $C$.
 Assume that the following algorithmic problems are decidable:

\begin{itemize}
  \item The Search Membership Problem for the subgroup $C$ in $A$ and in $B$.
  \item The Coset Representative Search Problem for the subgroup $C$ in $A$ and in $B$.
  \item The Conjugacy Membership Search Problem for $C$ in $A$ and in $B$.
  \item The Conjugacy Search Problem in $A$ and in $B$.
  \end{itemize}
   Then the Conjugacy Search Problem in $G$ is decidable for all pairs
$(g,h) \in G \times G$, where $g$ has cyclically reduced length 0.
\end{lemma}
 \begin{proof} Using Algorithm III (Theorem \ref{th:AlgII}) one can
 find cyclically reduced forms of  a given pair of elements $(g,h) \in G \times G$.
  In particular,
 one can verify if the cyclically reduced length of $g$ and $h$ is equal to 0.
  If so, then $g, h \in A \cup B \setminus C$. By the conjugacy criterion
  $g, h$ belong to one and the same factor $A$ or $B$, and they are conjugates there.
  Since the Search Conjugacy Problem is decidable in $A$ and $B$ the last condition is decidable.
  This proves the lemma.
 \end{proof}

\begin{remark}
\label{re:eff} The decision algorithms from Theorem \ref{th:main}
and Lemmas \ref{le:reg-0} and \ref{le:reg-1},  have polynomial time
complexity ``modulo" the algorithms for finding normal forms of
elements and the decision algorithms for the problems listed in the
statements.
\end{remark}

Combining Theorem \ref{th:main},  Lemmas \ref{le:reg-0} and
\ref{le:reg-1}, \ref{le:reg-C} and Corollary \ref{le:CR} altogether
 one can get the following general result.

\begin{theorem}
\label{th:main2}
 Let $G = A \ast_C B$ be a free product of finitely
presented groups $A$ and $B$ amalgamated over a finitely generated
subgroup $C$. Assume the following algorithmic problems are
decidable:

\begin{itemize}
  \item The  Membership Search Problem for $C$ in $A$ and in $B$.
  \item The Coset Representative Search Problem for the subgroup $C$in $A$ and $B$.
  \item The Cardinality Search Problem for $\Phi(Sub(C),A)$ in $A$ and for
   $\Phi(Sub(C),B)$ in $B$.
  \item The Conjugacy Search Problem in $A$ and in $B$.
  \item The Conjugacy Membership Search Problem for $C$ in $A$ and $B$.
   \item The Membership Problem for $N^\ast_A(C)$ and  $Z_A(C)$ in $A$,  and  $N^\ast_B(C)$
   and   $Z_B(C)$ in $B$.
\end{itemize}
Then the Conjugacy Search Problem in $G$ is decidable for arbitrary
pairs from $CR \times G$. \label{th:system}
\end{theorem}

\begin{corollary}
\label{co:freeCR} Let $G = A \ast_C B$ be  a free product  of free
groups $A$ and $B$ with amalgamated finitely generated subgroup $C$.
Then the Conjugacy Search Problem in $G$ is decidable for arbitrary
pairs from $CR \times G$.
\end{corollary}

\begin{corollary}
Let $G = A *_C B$ and $C$ is malnormal in $A$.  Assume the following
algorithmic problems are decidable:

\begin{itemize}
  \item The  Membership Search Problem for $C$ in $A$ and in $B$.
  \item The Coset Representative Search Problem for the subgroup $C$.
  \item The Cardinality Search Problem for $\Phi(Sub(C),A)$ in $A$ and for
   $\Phi(Sub(C),B)$ in $B$.
    \item The Conjugacy Membership Search Problem for $C$ in $A$ and $B$.
     \item The Conjugacy Search Problem decidable in $A$ and in $B$.
\end{itemize}
Then the Conjugacy Search Problem is decidable in $G$.
\end{corollary}

\begin{proof} Let $(g,h) \in G \times G$. Using Algorithm III
(Theorem \ref{th:AlgII}) one can
 find cyclically reduced forms of  the elements $g,h$. Assume for
 simplicity that $g$ and $h$ are cyclically reduced. If their
 cyclically reduced lengths are not equal, then they are not conjugates in $G$.
  Therefore we may assume that $l_0(g) = l_0(h)$.

1) Suppose  $l(g) = \l(h) \geq 2$. Since $C$ is malnormal in $A$
every element $g \in G$ with $l(g) \geqslant 2$ is regular  (see
Lemma  \ref{le:mal-deg}). Hence, in this case  by Theorem
\ref{th:main} the Conjugacy Search Problem for every pair $(g,h)$
with $l(g) \geqslant 2$  is decidable.

2) Suppose $l_0(g) = l_0(h) = 1$. In this case the argument from the
proof of Lemma \ref{le:reg-1} applies and gives the result.

3) Suppose $l_0(g) = l_0(h) = 0$, i.e., $g, h \in C$. By the
conjugacy criterion, there exists a sequence of elements $g = c_1,
c_2, \ldots, c_t = h$ from $C$ such that the neighboring elements
are conjugate either in $A$ or in $B$. By malnormality of $C$ in $A$
this implies that the neighboring elements are, in fact,  conjugate
in $B$. The latter is algorithmically decidable since the Conjugacy
Search Problem is decidable in $B$.

This proves the corollary.
\end{proof}

\end{document}